\documentclass[12pt,a4paper]{article}
\usepackage[english]{babel}
\usepackage{amssymb}
\usepackage{amsmath}
\usepackage[dvips]{graphicx}
\def\CC{\mathbb{C}}
\def\RR{\mathbb{R}}
\def\ZZ{\mathbb{Z}}

\def\A{{\mathcal A}}

\def\C{{\mathcal C}}

\def\L{{\mathcal L}}

\def\R{{\mathcal R}}
\def\K{{\mathcal K}}

\def\e{\varepsilon}
\def\O{{\mathcal O}}

\def\ov{\overline}
\def\a{\alpha}

\def\d{\delta}

\def\l{\lambda}
\def\LL{\Lambda}
\def\s{\sigma}
\def\S{\Sigma}
\def\De{\Delta}
 
\def\f{\varphi}
\def\wt{\widetilde}

\def\wh{\widehat}

\def\p{\partial}

\def\Mat{\operatorname{Mat}}
\def\Pf{\operatorname{Pf}}
\def\SL{\operatorname{SL}}
\def\GL{\operatorname{GL}}
\def\Sk{\operatorname{Sk}}
\def\Sym{\operatorname{Sym}}
\def\Sq{\operatorname{Sq}}
\def\Tr{\operatorname{Tr}}
\def\sq{\operatorname{sq}}
\def\sym{\operatorname{sym}}
\def\sk{\operatorname{sk}}
\def\sktr{\operatorname{sktr}}
\newtheorem{theorem}{Theorem}[section] 
\newtheorem{corollary}[theorem]{Corollary} 
\newtheorem{lemma}[theorem]{Lemma}

\newtheorem{example}[theorem]{Example} 
\newtheorem{prop}[theorem]{Proposition}
\newcommand{\bprop}{\begin{prop}}
\newcommand{\eprop}{\end{prop}}
\newtheorem{definition}[theorem]{Definition}
\newtheorem{rem}[theorem]{Remark}
\newcommand{\brem}{\begin{rem}}
\newcommand{\erem}{\end{rem}}
\newtheorem{rems}[theorem]{Remarks}
\newcommand{\brems}{\begin{rems}}
\newcommand{\erems}{\end{rems}}

\newtheorem{conj}[theorem]{Cojecture}
\newtheorem{problem}[theorem]{Problem}

\begin{document}

\title{{\bf Vanishing cycles of matrix singularities}}

\author{Victor Goryunov}
\date{}
\maketitle

\begin{abstract}  We study local singularities of holomorphic families of arbitrary square, symmetric and skew-symmetric matrices, that is, of mappings of smooth manifolds to the matrix spaces. Our main object is the vanishing topology of the pre-images of the hypersurface $\Delta$ of all degenerate matrices in assumption that the dimension of the source is at least the codimension of the singular locus of $\Delta$ in the ambient space. 

We start with showing that the complex link of $\Delta$ is homotopic to a sphere of the middle dimension and give a geometric interpretation of such spheres. This allows us to define vanishing cycles on the singular Milnor fibre of a matrix family, that is, on the local inverse image of $\Delta$ under a generic perturbation of the family. According to L\^e and Siersma such a fibre is a wedge of middle-dimensional spheres. We prove that in some important cases, which include the Damon-Pike conjecture and all simple matrix singularities, the number $\mu_\Delta$ of the spheres in the wedge 
is equal to the relevant Tjurina number $\tau$ of the family. 

We introduce two kinds of bifurcation diagrams for matrix families, and prove a Lyashko-Looijenga type theorem for the larger diagrams for all simple matrix families.

Making the first steps towards understanding the monodromy of matrix singularities, we define an intersection form on the singular Milnor fibres of corank 2 symmetric matrix families, which yields a complete description of the monodromy of such fibres. A modification of this approach reveals a quite unexpected relationship between certain Shephard-Todd groups, simple odd functions and simple corank 3 matrix singularities, those forming a sporadic part of the entire simple matrix classification obtained by Bruce and Bruce-Tari.  

We conclude with a general $\mu_\Delta=\tau$ conjecture for matrix singularities. 

\ 


\end{abstract}

\medskip 
This paper is about the vanishing topology of holomorphic families of arbitrary square, symmetric and skew-symmetric matrices, that is, of mappings of smooth manifolds to the corresponding matrix spaces. Complete classifications of simple families of the first two kinds were obtained more than 15 years ago by Bruce and Tari \cite{B,BT}, and a partial simple classification of skew-symmetric matrices was done by Haslinger about the same time \cite{H}. Topological and algebraic analysis of their results for a small number of parameters -- when a generically perturbed family avoids singularities of the set $\Delta$ of all degenerate matrices -- was carried out in \cite{BGZ,GZ,GM}. In particular, the monodromy of the Milnor fibre of the inverse image of the discriminant $\Delta$ was described in those papers, and a relation between the rank of the vanishing homology of such a fibre and the relevant Tjurina number of a matrix family was obtained. The last relation was actually proved in \cite{GM} not just for the matrix singularities but for sections of a wide class of singular hypersurfaces.    

Moving to a higher number of parameters in the matrix families brings a difference between two types of related Milnor fibres which have been diffeomorphic previously: the smooth Milnor fibre of a family which is the inverse image of the Milnor fibre of the discriminant, and the {\em singular Milnor fibre\/} which is the inverse image of $\Delta$ itself under a generic perturbation of the family. 

Like Milnor fibres of isolated function singularities, singular Milnor fibres of the matrix families are wedges of spheres of the middle dimension. This follows from \cite{Le,S}. However, being highly singular, they have been promising complicated calculations to yield any quantitative information, for example, the numbers $\mu_\Delta$ of the spheres in their wedges (called the {\em singular Milnor numbers}). 
One of the possible approaches to understand their topology was developed by Damon and Pike in \cite{DPI, DPII} where they introduced a rather involved machinery to calculate the Euler characteristic of such fibres via sections corresponding to fixed flags in the configuration spaces. The method is based on a choice of an appropriate sequence of free divisors and extends the L\^e-Greuel inductive procedure to calculate the Milnor number of an isolated complete intersection singularity \cite{LG}.  Nevertheless, numerous questions related to the geometry of bifurcation diagrams of the matrix families, to the monodromy of these families, to their relations with other types of singularities and so on have stayed untouched.

\medskip
The aim of this paper is to initiate a monodromy study of the singular Milnor fibres of the matrix families. 
The structure of the paper is as follows.

Section \ref{Seq} recalls the equivalences of matrix families we are going to work with.

In Section \ref{Slinks} we show that the complex link of the discriminant $\Delta$
in the space of all arbitrary square, symmetric or skew-symmetric matrices is homotopic to a sphere of the middle dimension, and describe these spheres in geometric terms. This allows us to define in Section \ref{Smvc} vanishing cycles on singular Milnor fibres of our matrix singularities, and prove in Section \ref{SDP} the Damon-Pike conjecture on the equality of the Tjurina and singular Milnor numbers for matrix families of a special type. Section \ref{Sbm} recalls the classification of simple matrix families from \cite{B,BT,H}. It also
extends the above $\mu_\Delta=\tau$ result to a wider class of matrix singularities, namely, to those corresponding to Arnold's functions on manifolds with boundaries. The latter, in particular, includes nearly all simple matrix singularities.

In Section \ref{Sgbd} we define two natural versions of bifurcation diagrams for matrix families, and prove a Lyashko-Looijenga type theorem for the majority of simple matrices.

Section \ref{S2x2} introduces a way to understand the monodromy of singular Milnor fibres which provides us with a complete description of such a monodromy for corank 2 symmetric matrices. Our approach may be considered as a generalisation of a modification of Arnold's construction of an intersection  form for boundary function singularities (cf. \cite{Ab}).

In Section \ref{Scork3}, our attention is on corank 3 simple symmetric matrix families, the only part of the simple classification not covered in earlier sections.
A higher-corank adjustment of the methods of the previous section
reveals a surprising relationship between these simple matrix singularities, odd functions and
certain Shephard-Todd groups.
In Section \ref{Scork3} we also notice that for these families the $\mu_\Delta = \tau$ equality holds as well, which extends the validity of the equality to all simple matrix singularities and also allows us to extend the Lyashko-Looijenga type theorem to all such singularities. 

\medskip
We conclude the paper with a general $\mu_\Delta = \tau$ conjecture for all three types of our matrix families.
Of course, another attractive theme emerging from the paper is a problem of description of the monodromy of singular Milnor fibres of matrix families to go beyond the results of Section \ref{S2x2}.

\medskip
{\em Acknowledgements.} I am grateful to David Mond and Oleg Karpenkov for very useful discussions. 
I am also grateful to the anonymous referees of the paper for their constructive suggestions which helped to improve the exposition.

\section{Matrix equivalences} \label{Seq}

A {\em matrix family\/} in this paper will be  a holomorphic mapping $M: U \to \Mat_n,$ where $U$ is an open domain in  $\CC^s,$ and $\Mat_n$ is one of three spaces of complex $n \times n$ matrices:
\begin{itemize}
\item[] $\Sq_n,$ arbitrary square matrices,
\item[] $\Sym_n,$ symmetric matrices, 
\item[] $\Sk_n,$ $n=2k$, skew-symmetric matrices.
\end{itemize}

In the $\Sq_n$ case, we say that two matrix families, $M_1$ and $M_2,$ are  {\em $\GL$-equivalent} if there exist a biholomorphism $\f$ of the source and two holomorphic maps $A,B: U \to \GL_n(\CC)$ such that
$$
M_1 \circ \f = A^\top M_2 B\,.
$$
For the $\GL$-equivalence  in the symmetric and skew-symmetric cases, we require existence of only one holomorphic $s$-parameter family $A$ of invertible $n \times n$ matrices such that
$$
M_1 \circ \f = A^\top M_2 A\,.
$$
Restriction to mappings $A,B: U \to \SL_n(\CC)$  in all three cases provides notions of the {\em $\SL$-equivalences}.

The relation between the matrix $\GL$- and $\SL$-equivalences is similar to that between the contact ($\K$) and right ($\R$) equivalences of holomorphic functions.

\bigskip
Germ versions of the above definitions are straightforward. 
A relevant equivalence class of germs of matrix families will be called a {\em matrix singularity}.  We denote by $\O_s$ the space of all holomorphic function germs on $(\CC^s,0),$ and identify the space of all holomorphic map germs $M:  (\CC^s,0) \to \Mat_n$ with the $\O_s$-module $\O_s^N$ where $N= {\rm dim}_\CC\Mat_n.$ 

Let $E_{j\ell}$ be the $n\times n$ matrix with the $j\ell$-entry 1, and all other entries zero. The extended tangent spaces to the $\GL$-equivalence classes of germs  $M:  (\CC^s,0) \to \Mat_n,$ with the source coordinates $x_1,\dots,x_s,$ are 
$$
\begin{array}{rcl}
T_{\GL,\Sq_n}M & = & \O_s \left<\p M/\p x_i, \   E_{j\ell}M, \  ME_{pq} \right>_ { i=1,\dots,s; \   j,\ell,p,q=1,\dots,n},  \qquad {\rm and}
\\ \\
T_{\GL,\Mat_n}M & = & \O_s \left<\p M/\p x_i, \   E_{j\ell}M + ME_{\ell j} \right>_{ i=1,\dots,s; \   j,\ell=1,\dots,n}
\quad
{\rm for} \ \ \Mat_n=\Sym_n,\Sk_n.
\end{array}
$$
The extended tangent spaces for the $\SL$-equivalences are contained in these: all the diagonal matrices $E_{jj}$ and $E_{pp}$ in the above expressions should be replaced, for example, by the differences $E_{jj}-E_{11}$ and $E_{pp}-E_{11}.$

We will denote by $\tau_{\GL,\Mat_n}(M)$ and $\tau_{\SL,\Mat_n}(M)$ the corresponding {\em Tjurina numbers\/} of $M,$ that is, the codimensions of the above extended tangent spaces in $\O_s^N.$

\medskip
There is an important situation when we have no difference between the two kinds of the extended tangent spaces, and hence between the two kinds of the Tjurina numbers. Namely,
we call a matrix family  $M:  (\CC^s,0) \to (\Mat_n,0)$ {\em quasi-homogeneous\/} if it is possible to assign positive weights $w_1, \dots, w_s$ to the source coordinates so that each entry $m_{ij}(x_1,\dots,x_s)$ would be a quasi-homogeneous function of degree $d_{ij}$ and each $2 \times 2$ minor would also be quasi-homogeneous, that is, $d_{ij} + d_{pq} = d_{iq} + d_{pj}$ must hold for the entire range of the indices. The last condition implies in particular that $\det(M)$ or $\Pf(M)$ are quasi-homogeneous too. All normal forms of matrix families appearing in this paper are quasi-homogeneous. 

\bprop \label{Pqh} 
The extended tangent spaces to the $\GL$- and $\SL$-equivalence classes of a quasi-homogeneous matrix family $M$ coincide. 
\eprop

{\em Proof.} 
Our quasi-homogeneity constraints on the degrees of the entries may be rewritten as
$ d_{ij} - d_{pj} = d_{iq} - d_{pq},$
that is, the difference between the degrees of the entries of two fixed rows is the same in all columns. Similar consistency holds when we compare entries in two fixed columns. 

In the $\Sq_n$ case, this observation allows us to split the degrees into sums of two sets of summands,  $d_{ij}=\d_i + \d'_j.$ So, we have for the Euler vector field $e = w_1 x_1 \p_{x_1} + \dots +  w_s x_s \p_{x_s}$: 
\begin{equation}
eM = \left((\d_i + \d'_j)m_{ij}\right)= \frak{d}M + M \frak{d}', \  {\rm with \  }   \frak{d} = diag(\d_1, \dots, \d_n), \  \frak{d}' = diag(\d'_1, \dots, \d'_n).
\label{Eqh}
\end{equation}
The ambiguity of shifting the whole $\d$-set by $d$ and the whole $\d'$-set by $-d$ corresponds to the kernel element $(I_n,-I_n)$ of the $\frak{gl}_n \times \frak{gl}_n$ action: \ \  $I_n M - MI_n = 0.$

The span of $\frak{sl}_n \times \frak{sl}_n$ and $(I_n,-I_n)$ has codimension 1 in $\frak{gl}_n \times \frak{gl}_n$: in terms of elements $(\frak{a},\frak{b})$ of the latter, this span is the hyperplane 
$trace(\frak{a}+\frak{b})=0.$ The element $(\frak{d},\frak{d}')$ is outside it. Relation (\ref{Eqh}) shows that its action is replaced by that of the Euler field to provide the effect of the action of the entire  $\frak{gl}_n \times \frak{gl}_n$ in the quasi-homogeneous $SL_n$ situation. 

\medskip
In the $\Sym_n$ and $\Sk_n$ cases, 
one similarly uses the symmetric splitting $d_{ij}=\d_i + \d_j.$  
 \hfill{$\Box$}

\bigskip
All our groups of matrix equivalences are in Damon's class of geometric subgroups of the $\K$-equi\-valence group (see \cite{D}), and therefore $\GL$- and $\SL$-miniversal deformations of a matrix family germ $M$ may be written as
$$
M + \l_1 \f_1 + \dots + \l_\tau \f_\tau\,,
$$ 
where $\tau$ is the relevant Tjurina number of $M,$ and the $\f_i \in \O_s^N$ form a basis of the quotient of $\O_s^N$ by the corresponding extended tangent space.

\bigskip 
By the {\em matrix corank\/} of a germ $M:  (\CC^s,0) \to \Mat_n$ we will understand the corank of the matrix $M(0).$ In  our symmetric and arbitrary square settings, a matrix corank $c$ family is equivalent to a family 
$\left( \begin{array}{cc}
M' & 0 \\ 0 & I_{n-c}
\end{array}
\right)$ where $M'$ is a germ of a $c \times c$ matrix family of the same type and $M'(0)$ is the zero matrix.
Similar reduction exists in the skew-symmetric case, with the only difference that the identity corner should be replaced by the block-diagonal matrix $J_{n-c}$ with the elementary  blocks
$J_2=\left( \begin{array}{cc} 0 & 1 \\ -1 & 0 \end{array} \right)$ along the diagonal.

Two germs of matrix families $M_i: (\CC^{s},0) \to \Mat_{n_i},$ $i=1,2,$  will be called {\em stably $\GL$-} or {\em $\SL$-equivalent} if there exists $n$ such that the two `extended'  families
$\left( \begin{array}{cc}
M_i & 0 \\ 0 & I_{n-n_i}
\end{array}
\right)$, or respectively $
\left( \begin{array}{cc}
M_i & 0 \\ 0 & J_{n-n_i}
\end{array}
\right)$, are equivalent.

\bigskip
In what follows we will need convenient skew-symmetric analogues of the standard trace and eigenvalues of a square matrix. So, we set
\begin{itemize}
\item the {\em skew trace} of $A \in \Sk_{2k}$ to be $\sktr(A) = \sum_{i=1}^k a_{2i-1,2i},$ and 
\item the {\em skew eigenvalues} of such an $A$ to be the solutions $\l_1, \dots, \l_k$ of the characteristic equation $\Pf (A-\l J_{2k})=0.$
\end{itemize}
We will also need a special kind of skew-symmetric matrices which we call {\em quaternionic\/}. These are matrices $A \in \Sk_{2k}$ in which all the $2 \times 2$ cells
$
\left(
\begin{array}{cc}
a_{2i-1,2j-1} & a_{2i-1,2j} \\
a_{2i,2j-1} & a_{2i,2j}
\end{array}
\right)
$
are additionally required to be of the form
$
\left(
\begin{array}{ccl}
z & w \\
-\ov w & \ov z
\end{array}
\right),
$  
$z,w \in \CC.$ In particular, all diagonal $2\times 2$ blocks 
$
\left(
\begin{array}{cc}
0 & a_{2i-1,2i} \\
- a_{2i-1,2i} & 0 
\end{array}
\right)
$
of a quaternionic matrix are real.

\begin{prop} \label{Pq}
All skew eigenvalues of a quaternionic matrix $A \in \Sk_{2k}$ are  real.
\end{prop}

{\em Proof.} For convenience, apply the same permutation to rows and columns of $A$ to bring it to the 4-block form
$$
\wh A = \left( \begin{array}{cc}
Z & W \\
-\overline W & \overline Z
\end{array}
\right) \in \Sk_{2k}\,,
$$
with $Z$ skew-symmetric and $W$ Hermitian. The same permutations transform $J_{2k}$ into
$\wh J_{2k} =  \left( \begin{array}{cc}
0 & I_k \\
-I_k & 0
\end{array}
\right). $
Therefore, skew eigenvalues of A are roots (all of them double) of the characteristic equation
$$
\det(\wh A - \l \wh J_{2k})=0, \ \ \  {\rm that \  is,} \ \ \  \det(\wh A \wh J_{2k}^{-1} - \l I_{2k})=0\,,
$$
where $\wh A \wh J_{2k}^{-1}= -\wh A \wh J_{2k} =
-\left( \begin{array}{cc}
W & -Z \\
\overline Z & \overline W
\end{array}
\right)
$ is Hermitian due to our constraints on $Z$ and $W.$
 
\hfill{$\Box$}

\section{Complex links of matrix discriminants} \label{Slinks}
We denote by $\De \subset \Mat_n$ the set of all degenerate matrices, and call this set the {\em discriminant}.
Under the {\em discriminant $\De (M)$ of a particular matrix family $M$} we understand the zero set of the function $\det \circ M$ or $\Pf \circ M,$ 
that is, the inverse image of $\De$ under the mapping $M.$

We use the notation $\Tr_\a$ for the hyperplane in $\Mat_n$ of all matrices with the (skew) trace $\a.$ We write $\De^{\sym},$ $\Tr^{\sq}_\a,$ etc. if we want to emphasise which of the three cases we are considering.

\medskip
According to \cite{B,BT,H},  a generic map germ $(\CC^{N-1},0) \to (\Mat_n,0)$ is --- up to any of our equivalences --- an embedding whose image is the hyperplane germ $(\Tr_0,0).$
Moreover, $\GL$- and $\SL$-miniversal deformations of such maps  are one-parameter families of embeddings of $\CC^{N-1}$ into $\Mat_n$ as the hyperplanes $\Tr_\a,$
$\a \in (\CC,0).$ Therefore, the complex links of the matrix discriminants are the sets of all degenerate matrices with a fixed non-zero (skew) trace.

\begin{theorem} \label{Tlink}
The complex link of the discriminant $\De \subset \Mat_n$ is homotopy equivalent to a sphere $S^{N-2}.$ If $\l$ is a positive real number, then this sphere may be taken to be the set of  respectively
\begin{itemize}
\item[(Sym)] all degenerate real matrices $A \in \Sym_n$ with the trace $n\l$ and all eigenvalues non-negative;
\item[(Sq)]  all degenerate hermitian matrices $A \in \Sq_n$ with the trace $n\l$ and all eigenvalues non-negative;
\item[(Sk)] all degenerate quaternionic matrices $A \in \Sk_n,$ $n=2k,$ with the skew trace $k\l$ and all skew eigenvalues non-negative.
\end{itemize}
\end{theorem}

{\em  Proof.} 
We consider the symmetric case in detail, and point out the adjustments required in the two other cases.

\medskip
{\em Symmetric matrices.} The only level with non-isolated critical points of the restriction of the determinantal function to the hyperplane $\Tr_{n\l}^{\sym},$ $\l \ne 0,$ is its zero level.
If we show that this restriction of the determinantal function has exactly one Morse critical point outside the zero level then the claim about the homotopy
type of the link will follow from \cite{Le,S}.

\begin{lemma} \label{Llink}
For $\l \ne 0,$ the only invertible matrix which is a critical point of the restriction of the determinantal
function to $\Tr_{n\l}^{\sym}$  is $\l I_{n}.$
\end{lemma}

{\em Proof of the lemma.}
We parametrise the hyperplane $\Tr_{n\l}^{\sym} \subset \Sym_n$ as a shift of $\Tr_0^{\sym}$ on which all the entries $x_{ij}$, $1 \le i \le j \le n,$  except for $x_{11},$ are taken as independent coordinates, that is, by the map $M_\l: \CC^{N-1} \to \Sym_n,$ $ N=n(n+1)/2:$
\begin{equation}
\label{Elink}
M_\l (x) = \l I_n + L^{\sym}(x), {\rm \ \,  where \  \,  }  
L^{\sym}(x)=\left(
\begin{array}{cccccccccccc}
 - \sum_{i=2}^{n} x_{ii}  & x_{12} & x_{13} & \dots & x_{1n} \\
 x_{12} & x_{22} & x_{23} & \dots & x_{2n} \\
x_{13} & x_{23} &  x_{33} & \dots & x_{3n} \\
 \vdots & \vdots & \vdots & \ddots & \vdots \\
x_{1n} & x_{2n} & x_{3n} & \quad\dots\quad & x_{nn} 
\end{array}
\right).
\end{equation}

Here the partial derivative $\p\det(M_\l(x)) / \p x_{ij},$ $i<j,$ is twice the $ij$-cofactor of the symmetric matrix $M_\l (x).$  Therefore, vanishing of these derivatives for all $i<j$ at an invertible matrix guarantees that  the inverse of such a matrix is diagonal. However, the block-diagonal structures of a matrix and of its inverse coincide.
Hence an invertible critical matrix itself is also diagonal: $x_{ij}=0$ for all $i<j.$  
This reduces the determinant of a member of the $M_\l$ family to the product of its diagonal entries subject to the  constraints that their sum $n\l$ is fixed and that all of them are non-zero. Such a product has only one critical point, when all its factors are equal, that is, all the $x_{ii}$ are zero.  
 \hfill{$\Box$} 

\medskip
We now look at the terms in $\det(M_\l(x))$ which are quadratic in $x.$ They are $\l^{n-2}$ times the sum of all $2 \times 2$ principal minors of $L^{\sym}.$ This sum reduces to 
\begin{equation}
\label{Equad}
-\sum\limits_{1\le i < j \le n} x_{ij}^2 -  \sum_{i=2}^n x_{ii}^2 - \sum_{2\le i < j \le n} x_{ii}x_{jj}\,.
\end{equation}
This quadratic form is non-degenerate, and therefore the critical point $M_\l(0)=\l I_n$ is indeed a Morse critical point of the determinantal function.

For $\l$ real and positive, the form (\ref{Equad}) is negative definite on the real part of $\Tr^{\sym}_{n\l}.$ This implies the {\em (Sym)} claim of the Theorem. Indeed, taking small $\e>0,$ we see that the  Morse vanishing sphere $S^{N-2}$ on the under-critical level $\det =\l^n - \e$ is real (see, for example, Section 1.3 in \cite{AGV2}), and it is sent by the negative gradient flow of the determinantal function on the real part of $\Tr_{n\l}^{\sym}$ to the set of real symmetric matrices specified in the Theorem.

\medskip
{\em Square matrices.} We parametrise $\Tr_{n\l}^{\sq}$ similar to (\ref{Elink}) with the only difference that the over- and under-diagonal entries are now independent:
\begin{equation}
\label{Elinksq}
M_\l(x)=\l I_n + L^{\sq}(x), {\rm \ \,  where \  \,  }  
L^{\sq}(x)=\left(
\begin{array}{cccccccccccc}
 - \sum_{i=2}^{n} x_{ii}  & x_{12} & x_{13} & \dots & x_{1n} \\
 x_{21} & x_{22} & x_{23} & \dots & x_{2n} \\
x_{31} & x_{32} &  x_{33} & \dots & x_{3n} \\
 \vdots & \vdots & \vdots & \ddots & \vdots \\
x_{n1} & x_{n2} & x_{n3} & \quad\dots\quad & x_{nn} 
\end{array}
\right), 
\end{equation}
$ x \in \CC^{n^2-1},$ $\l\ne 0.$
This replaces the first term in (\ref{Equad}) with $-\sum_{1\le i < j \le n} x_{ij} x_{ji}.$ After the replacement we set $x_{ji} = \bar x_{ij}$ for all $i \le j$ to treat the amended form as a negative definite quadratic form on a real space.

\medskip
{\em Skew-symmetric matrices.} For a parametrisation of $\Tr_{k\l}^{\sk}$ we take
\begin{equation}
\label{Elinksk}
M_\l (x)=\l J_{2k} + L^{\sk}(x),
\end{equation}
where
$$
L^{\sk}(x)=\left(
\begin{array}{cccccccccccc}
0 &\!\!\! - \sum_{i=2}^{k} x_{2i-1,2i}  \!\!\! & x_{1,3} & x_{1,4} & \dots &  x_{1,2k-1} & x_{1,2k} \\ 
 \sum_{i=2}^{n} x_{2i-1,2i}\!\!\!  &  0  &   x_{2,3} & x_{2,4} & \dots &  x_{2,2k-1} & x_{2,2k} \\ 
-x_{1,3} & -x_{2,3} & 0 &  \!\!\!   x_{3,4}  \!\!\! & \dots &  x_{3,2k-1} & x_{3,2k} \\ 
-x_{1,4} & -x_{2,4} & \!\!\!  -  x_{3,4} \!\!\! & 0 & \dots &  x_{4,2k-1} & x_{4,2k} \\ 
 \vdots & \vdots & \vdots & \vdots & \ddots & \vdots & \vdots \\ 
-x_{1,2k-1} & -x_{2,2k-1} & -x_{3,2k-1} & -x_{4,2k-1} & \dots &  0 &   x_{2k-1,2k} \\  
-x_{1,2k} & -x_{2,2k} & -x_{3,2k} & -x_{4,2k} & \,\dots\, &   -  x_{2k-1,2k}  & 0 
\end{array}
\right), 
$$
$ x \in \CC^{N-1},$ $ N=2k(2k-1)/2,$ $\l\ne 0.$

Similar to the proof of Lemma \ref{Llink}, vanishing of the derivatives of $\det(M_\l(x))= \linebreak \Pf^2(M_\l(x))$ with respect to all the $x$-coordinates except for the $x_{2i-1,2i},$ $i>1,$ implies that a critical invertible matrix must be block-diagonal, with the multiples of the $J_2$ blocks along the diagonal, that is,
its Pfaffian is $(\l - \sum_{i=2}^{k}  x_{2i-1,2i} )\prod_{i=2}^k(\l +  x_{2i-1,2i}).$ The only critical point of this product outside its zero level is $M_\l (0)=\l J_{2k}.$

The terms of $\Pf(M_\l (x))$ quadratic in $x$ are $\l^{k-2}$ times
$$
-\sum_{1\le i < j \le k} \left| 
\begin{array}{cc}
x_{2i-1,2j-1} & x_{2i-1,2j} \\
x_{2i,2j-1} & x_{2i,2j}
\end{array}
\right|
 -  \sum_{i=2}^k x_{2i-1,2i}^2 - \sum_{2\le i < j \le k} x_{2i-1,2i}x_{2j-1,2j}\,.
$$ 
The determinants here become positive definite quadratic forms on $\RR^4$ when we set  $x_{2i,2j} =  \bar x_{2i-1,2j-1}$ and  $x_{2i,2j-1} = - \bar  x_{2i-1,2j}.$ Taking also all the $x_{2i-1,2i}$ real, we end up with our skew-symmetric matrices being quaternionic. 
\hfill{$\Box$}

\brem\label{Rblock} {\em 
The observation used in the proof of Lemma \ref{Llink} about the coincidence of the block-diagonal structures of an invertible matrix $M$ and of its inverse will be also applied in our later calculations, when it is known that the cofactors of an appropriate subset of entries of $M$ are zero. 
}
\erem

\section{Vanishing cycles} \label{Svc}
From now on we assume that the dimension $s$ of the domains of our matrix families $\CC^s \to \Mat_n$ is such that {\em in general\/} their images meet the singular loci of the discriminants $\De \subset \Mat_n,$ that is, $s >2,3,5$ in respectively  symmetric, arbitrary square and skew-symmetric cases.

Let $M: (\CC^s,0) \to \Mat_n$ be a germ of a matrix family with a finite $\SL$ Tjurina number $\tau,$ and $\{M_\l\}_{\l \in \LL},$ $M_0=M,$ a representative of its $\SL$-miniversal deformation. According to \cite{Le,S},  for sufficiently small radii $\e>0$ and $\eta>0,$ there exist closed $2s$- and $2\tau$-dimensional balls $B_\e \subset \CC^s$ and $B_\eta \subset \LL,$ centred at the origins, such that the set $\De(M_\l) \cap B_\e$ is homotopic to a wedge of a finite number of $(s-1)$-dimensional spheres if $\l \in B_\eta.$

The variety $\De(M_\l) \cap B_\e$ is the zero level of the function $\det \circ M_\l$ or $\Pf \circ M_\l$ on $B_\e.$ Due to the constraints on the dimension $s,$ this is the only level of the function which is critical for all $\l \in \LL$ if the corank of $M(0)$ is higher than the lowest positive, that is, greater than 1 in the $\Sq$ and $\Sym$ cases and greater than 2 in the $\Sk$ case. Moreover, the singular locus of the zero level is of positive dimension  if $s$ is not the minimal we have allowed.
Critical points $p_i,$
$i=1,\dots,r,$ on all other levels are isolated, and therefore the function has the corresponding ordinary Milnor numbers $\mu_i$ at them. By \cite{Le,S}, the number of the spheres in the wedge of the set $\De(M_\l) \cap B_\e$ is $\mu_1 + \dots + \mu_r \,.$

 For a generic $\l \in B_\eta,$ the number of the spheres in the wedge achieves its maximum. 

\begin{definition} {\em 
We denote this maximum $\mu_\De (M)$ and call it {\em the singular Milnor number of the matrix family $M.$} The set $\De(M_\l) \cap B_\e$ homotopy equivalent to a wedge of $\mu_\De (M)$ copies of $S^{s-1}$ will be called {\em the singular Milnor fibre of the family $M.$} }
\end{definition}

For convenience, we will use the same terminology even if the corank of $M(0)$ is the lowest positive. Of course, in such cases the above $\De(M_\l) \cap B_\e$ is just the ordinary smooth Milnor fibre of the function $\det \circ M_\l$ or $\Pf \circ M_\l,$ and $\mu_\De (M)$ is its ordinary Milnor number.

\medskip
In this section we calculate singular Milnor numbers of matrix families related to isolated singularities of functions on manifolds with and without boundaries. 
In particular, they include nearly all simple matrix singularities from \cite{B,BT}. The only remaining simple singularities, of 4-parameter corank 3 symmetric matrices and of their 7-parameter square versions, will wait till the final section of this paper.

Later on, to avoid reminders about the choice of the balls, we will refer to the sets $\De(M_\l) \cap B_\e$ with $\l \in  B_\eta$ as {\em local\/} sets $\De(M_\l).$

\subsection{Matrix vanishing cycles} \label{Smvc}
In what follows, it will be more convenient for us to write miniversal deformations of the parametrisations $L: (\CC^{N-1},0) \to (\Mat_n,0)$ of the hyperplanes $\Tr_0$ used in (\ref{Elink}), (\ref{Elinksq}), (\ref{Elinksk})
in slightly different ways than in those formulas, namely as respectively
\begin{equation}
\label{E1}
L^{\sym}(x) + \l E_{11}, \quad L^{\sq}(x) + \l E_{11}, \quad L^{\sk}(x) + \l (E_{12}-E_{21})\,.
\end{equation}

As it has already been said, the parametrisations $L$
are the local $\GL$- and $\SL$-normal forms of generic hypersurface embeddings. 
So, they are also normal forms of the most generic rank $N-1$ linear parts of map germs  $M: (\CC^{s},0) \to (\Mat_n,0) \simeq (\CC^{N},0),$ $s=N-1+m.$ Hence any map $M$ with such a linear part  may be considered as an $m$-parameter deformation of the relevant family $L,$ and therefore it  reduces to the form (\ref{E1}) with $\l = g (z),$ where $z = (z_1,\dots,z_m)$ are the additional variables. 

In the simplest situation $g$ has a Morse critical point at the origin, and the matrix families are $\SL$- and $\GL$-equivalent to 
\begin{equation}
\label{Evs0}
\begin{array}{rcl}
L^{\sym}(x) & - & (z_1^2 + \dots + z_m^2) E_{11}, \\
L^{\sq}(x) & - & (z_1^2 + \dots + z_m^2) E_{11}, \\
L^{\sk}(x) & - & (z_1^2 + \dots + z_m^2) (E_{12}-E_{21})\,.
\end{array}
\end{equation}
For these families $\tau_{\SL,\Mat_n} = \tau_{\GL,\Mat_n} = 1,$ and we have the miniversal deformations
\begin{equation}
\label{Evs}
\begin{array}{rcl}
L^{\sym}(x) & + & (\l - z_1^2 - \dots -z_m^2) E_{11}, \\
L^{\sq}(x) & + & (\l - z_1^2 - \dots - z_m^2) E_{11}, \\
L^{\sk}(x) & + & (\l - z_1^2 - \dots -z_m^2) (E_{12}-E_{21})\,,
\end{array}
\end{equation}
where we are back to using $\l \in \CC$ as a deformation parameter. These are versal deformations of the only codimension 1 matrix singularities. Adjusting our proof of Theorem \ref{Tlink} to these deformations, we obtain 
the following description of all possible vanishing cycles appearing in singular Milnor fibres of our matrix singularities.

\begin{corollary}
The singular Milnor fibre of each of the matrix families in (\ref{Evs0}) is homotopy equivalent to a sphere $S^{s-1},$ $s=N-1+m,$ $m \ge 0.$
 If $\l$ is a fixed positive real number, then this sphere may be taken to be the set of  respectively
\begin{itemize}
\item[(Sym)] all degenerate real matrices $A \in \Sym_n$ with the trace $\l - z_1^2 - \dots -z_m^2,$ $z \in \RR^m,$ and all eigenvalues non-negative;
\item[(Sq)]  all degenerate hermitian matrices in $A \in \Sq_n$ with the trace $\l - z_1^2 - \dots -z_m^2,$ $z \in \RR^m,$ and all eigenvalues non-negative;
\item[(Sk)] all degenerate quaternionic matrices $A \in \Sk_n,$ $n=2k,$ with the skew trace $\l - z_1^2 - \dots -z_m^2,$ $z \in \RR^m,$ and all skew eigenvalues non-negative.
\end{itemize}
\end{corollary}

For a proof, we just need to notice that, for a fixed $\l \ne 0$, all families in (\ref{Evs}) contain only one critical point of the determinantal function off its zero level.
Indeed, for example in the symmetric case, the argument used in Lemma \ref{Llink} shows that a critical matrix is diagonal. After that the vanishing of the derivative of the determinant with respect to any $z_i$ yields $z_i=0$ outside $\Delta.$

Thus, all the vanishing cycles are the order $m$ suspensions of the complex links of the matrix discriminants.

\begin{definition} \label{Dcork} {\em
Vanishing cycles corresponding to the stable $\GL$- or $\SL$-equivalence classes of the families (\ref{Evs0}) of the matrix corank $n$ will be called} corank $n$ vanishing cycles. 
\end{definition}

In particular, corank 1 vanishing cycles in the symmetric and square cases, as well as corank 2 vanishing cycles in the skew-symmetric case are the ordinary Morse vanishing cycles of isolated hypersurface singularities.  Vanishing cycles of higher coranks have singularities, they are homeomorphic but not diffeomorphic to the spheres.

\brem {\em
We would like to emphasise that the `vanishing'  in the last definition takes place within the one-parameter sets of the {\em singular} Milnor fibres corresponding to the local deformations (\ref{Evs}) when $\l$ tends to zero. }
\erem

\subsection{Damon-Pike conjecture} \label{SDP}
We will now prove

\begin{theorem}\label{TDP}
Consider matrix singularities $M: (\CC^{s},0) \to (\Mat_n,0) \simeq (\CC^{N},0)$ with a generic rank $N-1$ linear part, that is,
\begin{equation}
\label{E2}
L^{\sym}(x) + g(z) E_{11}, \quad L^{\sq}(x) + g(z) E_{11}, \quad L^{\sk}(x) + g(z) (E_{12}-E_{21})\,,
\end{equation}
$ x \in \CC^{N-1},$ $g \in \O_m.$ Assume the $\SL$ Tjurina number of such a matrix family is finite.  Then
$$\mu_\Delta(M) = \tau_{\SL,\Mat_n} (M) = \mu(g)\,.$$
\end{theorem}

The $\mu(g)$ here is the standard Milnor number of the isolated function singularity $g.$

\medskip
This property was conjectured by Damon and Pike in \cite{DPII}, with a quasi-homogeneity of $g$ requested for the last equality.

\ 

{\em Proof.}  
We consider only the symmetric case, the other two being its word-to-word repetitions.

So, take a small generic deformation $\{M_t\}_{t \in (\CC,0)}$ of 
$$M_0(x,z) = M(x,z) = L^{\sym}(x) + g(z) E_{11}.$$ 
Following the approach used in the previous subsection, we can set 
$$M_t (x,z) = L^{\sym}(x) + g_t(z) E_{11},$$
 where $\{g_t\}_{t \in (\CC,0)}$ is a small generic deformation of $g_0=g.$ 
To prove the equality of $\mu_\Delta(M)$ and $\mu(g),$ we need to show that, for a fixed $t \ne 0,$ the function $\det \circ M_t$ has exactly $\mu(g)$ Morse critical points close to the origin and outside its zero level $\De(M_t).$ 

Following the argument used in Lemma \ref{Llink} and highlighted in Remark \ref{Rblock}, we see that a critical matrix in the $M_t$ family off $\Delta$  must be diagonal and hence has determinant
$$
\left( -\sum_{i=2}^n x_{ii} + g_t(z) \right) \cdot \prod_{i=2}^n x_{ii}\,.
$$ 
Vanishing of the derivatives of this expression with respect to all the $x_{ii}$ outside zeros of the expression itself implies
$$
x_{22} = x_{33} = \dots = x_{nn} = g_t(z)/n\,.
$$
Therefore the determinant of a critical matrix is $\left(g_t(z)/n\right)^n.$ Critical points of the last function outside its zeros are exactly critical points of $g_t.$ Due to the genericity of the deformation $\{g_t\}$ of $g_0=g,$ these points are all Morse and their number is $\mu (g).$ 

\medskip
The remaining part of the theorem, that $\tau_{\SL,\Sym_n} (M) = \mu(g),$ follows from a calculation exercise (cf. \cite{B}) showing that
$$
\O_{N-1+m}^N / T_{\SL,\Sym_n}M \simeq \left( \O_m / \O_m \left< \p g/\p z_1,\dots,  \p g/\p z_m \right> \right) E_{11}\,.
$$

\vspace{-25pt}\hfill{$\Box$}

\ 

The last isomorphism holds also for the square matrices in (\ref{E2}), and --- with replacement of $E_{11}$ by $(E_{12} - E_{21})$ --- in the skew-symmetric case too. Due to that, $\SL$-miniversal deformations of all matrix families (\ref{E2}) may be obtained by replacing the function $g$ in them by its $\R$-miniversal deformation.  

In the $\GL$ situation, for the matrix families (\ref{E2}) we have 
$ \tau_{\GL,\Mat_n} (M) = \tau(g),$  where $\tau(g)$ is the Tjurina number for the $\K$-equivalence of functions (see Proposition 4.9 in \cite{B}). Respectively, for $\GL$-miniversal deformations of these families, we replace the function $g$ by its $\K$-miniversal deformation. 

\begin{rem} {\em
The reduction of the $\SL$- and $\GL$-classifications of the matrix families (\ref{E2}) to respectively the $\R$- and $\K$-classifications of the function germs $g: (\CC^s,0) \to (\CC,0)$ is in fact what is very much expected. Indeed, according to what was said at the beginning of the previous subsection, the $\SL$-equivalences leave the whole group of diffeomorphisms of $(\CC^s,0)$ for normalising $g.$ In the $\GL$ cases we additionally have an opportunity of multiplying the whole matrix by functional multiples $h(z) I_n$ of the identity matrix, $h(0) \ne 0.$ Modulo the same rescaling of all the $x$-coordinates, the last move is just multiplication of $g$ by $h$. The hypersurface $g=0$ in $\CC^m$ has a clear geometric meaning: it is the inverse image of the zero matrix under the mapping $M.$ }
\end{rem}

\subsection{$\mu = \tau$ theorem for a matrix version \\ of boundary function singularities} \label{Sbm}

The type of matrix families we are studying in this subsection generalises that from the Damon-Pike conjecture, and contains nearly all simple matrix classes, at least in the symmetric and square settings. The families we are considering now have a more general top left corner:
\begin{equation}
\label{Esymb}
M^{\sym}(x,z) = \left(
\begin{array}{cccccccccccc}
 - \sum_{i=3}^{n} x_{ii} & x_{12} & x_{13} & \dots & x_{1n} \\
 x_{12} & x_{22} & x_{23} & \dots & x_{2n} \\
x_{13} & x_{23} &  x_{33} & \dots & x_{3n} \\
 \vdots & \vdots & \vdots & \ddots & \vdots \\
x_{1n} & x_{2n} & x_{3n} & \quad\dots\quad & x_{nn} 
\end{array}
\right) + h(x_{22},z)E_{11}\,, 
\end{equation}
\begin{equation}
\label{Esqb}
M^{\sq}(x,z)=\left(
\begin{array}{cccccccccccc}
 - \sum_{i=3}^{n} x_{ii} & x_{12} & x_{13} & \dots & x_{1n} \\
 x_{21} & x_{22} & x_{23} & \dots & x_{2n} \\
x_{31} & x_{32} &  x_{33} & \dots & x_{3n} \\
 \vdots & \vdots & \vdots & \ddots & \vdots \\
x_{n1} & x_{n2} & x_{n3} & \quad\dots\quad & x_{nn} 
\end{array}
\right) + h(x_{22},z)E_{11}\,, 
\end{equation}
\begin{multline}
\label{Eskb}
M^{\sk}(x,z)= \left(
\begin{array}{cccccccccccc}
0 &\!\!\! - \sum_{i=3}^{k} x_{2i-1,2i}  \!\!\! & x_{1,3} & x_{1,4} & \dots &  x_{1,2k-1} & x_{1,2k} \\ 
 \sum_{i=3}^{n} x_{2i-1,2i}  \!\!\!  &  0  &   x_{2,3} & x_{2,4} & \dots &  x_{2,2k-1} & x_{2,2k} \\ 
-x_{1,3} & -x_{2,3} & 0 &  \!\!\!   x_{3,4}  \!\!\! & \dots &  x_{3,2k-1} & x_{3,2k} \\ 
-x_{1,4} & -x_{2,4} & \!\!\!  -  x_{3,4} \!\!\! & 0 & \dots &  x_{4,2k-1} & x_{4,2k} \\ 
 \vdots & \vdots & \vdots & \vdots & \ddots & \vdots & \vdots \\ 
-x_{1,2k-1} & -x_{2,2k-1} & -x_{3,2k-1} & -x_{4,2k-1} & \dots &  0 &   x_{2k-1,2k} \\  
-x_{1,2k} & -x_{2,2k} & -x_{3,2k} & -x_{4,2k} & \,\dots\, &   -  x_{2k-1,2k}  & 0 
\end{array}
\right) \\ \\ + h(x_{3,4},z)(E_{12}-E_{21})\,.
\end{multline}
Here again $x \in \CC^{N-1}$ and  $z \in \CC^m,$ while $h \in \O_{m+1}.$ Note that the summations in the matrices are one term shorter now compared with the matrices $L$ from the previous subsection. Of course, we are back to the situation of that subsection if the derivative of $h$ at the origin by its only $x$ argument is not zero. 

\medskip
We now recall the $\GL$-simple classification under our earlier assumption that the number $s=N-1+m$ of the matrix parameters should be high enough for the singular Milnor fibre to be indeed singular. The main part of this classification is common to all the sizes, and is related to the classification of the germs $h \in \O_{m+1}$ in (\ref{Esymb}--\ref{Eskb})
up to the $\R_\p$-equivalence of functions on the manifold $\CC^{m+1}$ with the boundary $\CC^m$ given by the vanishing of the only $x$ coordinate in $h$ \cite{Ab}. 

\begin{theorem}\label{Tsimple}
\begin{itemize}
\item[(a)] 
\cite{B,BT,H} $\GL$-simple matrix singularities of the forms (\ref{Esymb}), (\ref{Esqb}) and (\ref{Eskb}) are classified by the $\R_\p$-simple types $X_\tau = A_\tau,D_\tau,E_\tau,B_\tau,C_\tau,
F_4$  of the function germs $h.$

\item[(b)]  
\cite{B} Up to the stable $\GL$-equivalence, the only other simple classes in the symmetric case are the following corank 3 matrix families in four variables: 
$$
\begin{array}{l|l|l|l}
I_{k+1}, \  k \ge 1 & II_4 & II_5 & II_6 \\
\hline \vspace{-10pt}
&&& \\   
\left( 
\begin{array}{ccc}
x & 0 & z \\
0 & y + x^k & w \\
z & w & y
\end{array}
\right)  
 &  
\left( \begin{array}{ccc}
x & w^2 & y \\
w^2 & y & z \\
y & z & w 
\end{array}
\right)  
&  
\left( \begin{array}{ccc}
x & 0 & y\!+\!w^2 \\
0 & y & z \\
y\!+\!w^2 & z & w 
\end{array}
\right)   &  
\left( \begin{array}{ccc}
x & w^3 & y \\
w^3 & y & z \\
y & z & w 
\end{array}
\right)  
\end{array}
$$

\item[(c)] 
(\cite{BT}, slightly corrected)  Up to the stable $\GL$-equivalence, the only other simple classes in the square case are corank 3 families in seven variables obtained from the symmetric matrices in item (b) by addition of the generic skew-symmetric matrix family \linebreak  
$
U=
\left( \begin{array}{ccc}
0 & u_{12} & u_{13} \\
-u_{12} & 0 & u_{23}  \\
-u_{13} & -u_{23} & 0
\end{array}
\right)
$
in three new variables.  
\end{itemize}
\end{theorem}

The families of item (a) will be denoted $X_\tau^{\Mat_n},$ for all our three types of matrices. Their $\GL$ Tjurina numbers are $\tau.$ For example, the families in (\ref{Evs0}) are the $A_1^{\Mat_n}$ singularities, and their one-parameter deformations (\ref{Evs}) are their $\GL$-miniversal deformations. The singularities of item (b) are also indexed by their $\GL$ Tjurina numbers, and these numbers are the same for the corresponding seven-variable families of square matrices. 
It is not known if the skew-symmetric simple classification goes beyond the part of item (a) of the form (\ref{Eskb}).
The $\SL$-simple classification coincides with the $\GL$ one (including the last incompleteness comment) due to the quasi-homogeneity of all matrix singularities of the theorem. Of course, $\GL$-non-simplicity implies $\SL$-non-simplicity since the $\SL$-classification is finer.

\medskip
Thus, matrix singularities of the forms (\ref{Esymb}--\ref{Eskb}) are playing a special role in the classification, and our aim now is to extend the Damon-Pike conjecture to such families.

So, let us find the singular Milnor numbers for matrices  (\ref{Esymb}--\ref{Eskb}).
According to \cite{B,BT} and similar skew-symmetric calculations, all possible deformations of these matrix families may be obtained by deforming the functions $h.$ 
To make some difference with the previous subsection, we shall now consider the skew-symmetric case.
Therefore, let $\{h_t\}_{t \in (\CC,0)}$ be a generic small deformation of $h_0=h$ in (\ref{Eskb}) and $\{M_t\}_{t \in (\CC,0)}$ the corresponding matrix deformation. Using Remark \ref{Rblock}, we see that
vanishing of the derivatives of the $\Pf \circ M_t$ outside $\Delta(M_t)$ with respect to all the $x_{j,\ell}$ except for the $x_{2i-1,2i},$ $i>1,$ implies that a non-degenerate critical matrix in the family  $M_t,$ $t\ne 0,$ has the same block-diagonal structure as the standard matrix $J_{2k}$ and therefore its Pfaffian  
is the product 
$$  
\left( - \sum_{i=3}^{k} x_{2i-1,2i} + h_t(x_{3,4},z)\right) \prod_{i=2}^k x_{2i-1,2i}\,.
$$
Vanishing of  the derivatives of this product with respect to all the $x_{2i-1,2i},$ $i>2,$ implies that the values of all these coordinates must be the same and equal to $h_t(x_{3,4},z)/(k-1).$ This leaves us with the Pfaffian equal to
\begin{equation}
\label{Ecrbound}
\left( h_t(x_{3,4},z)/(k-1) \right)^{k-1} \cdot x_{3,4}\,.
\end{equation}
Vanishing of the derivatives with respect to all the still remaining coordinates gives us the final conditions
$$
\p h_t/ \p z_1 = \p h_t/ \p z_2 = \dots = \p h_t/ \p z_m = h_t + (k-1) x_{3,4} \p h_t/ \p x_{3,4} = 0 \,.
$$ 
Therefore,  if the quotient 
\begin{equation}
\label{Ebound}
Q_h =
\O_{m+1}/ \O_{m+1} \left< {\p h \over\p z_1}, \dots,  {\p h \over\p z_m}, h + (k-1) x_{3,4} {\p h \over\p x_{3,4}} \right>
\end{equation}
is finite-dimensional, then  $\mu_\Delta (M) = {\rm dim}_\CC \,Q_h.$

The expression for the local algebra $Q_h$  in the symmetric and square cases has $n$ instead of $k,$ and $x_{22}$ instead of $x_{3,4}.$

On the other hand, direct calculations (cf. \cite{B,BT}) of the extended tangent spaces to the $\SL$-equiva\-lence classes of the families (\ref{Esymb}),  (\ref{Esqb}) and (\ref{Eskb}) show that 
$\O_{N-1+m}^N / T_{\SL,\Mat_n}M$ is isomorphic to either $Q_h  E_{11}$ or respectively $Q_h  (E_{12}-E_{21}).$

Thus we have proved

\begin{theorem}\label{Tb}
Assume the $\SL$ Tjurina number of a matrix family $M$ of a form (\ref{Esymb}), (\ref{Esqb}) or (\ref{Eskb}) is finite. Then

\vspace{-5pt}
\centerline{$\mu_\Delta(M) = \tau_{\SL,\Mat_n} (M) =  {\rm dim}_\CC \,Q_h \,.$}
\end{theorem} 

\begin{rem} {\em 
If $h$ is quasi-homogeneous then ${\rm dim}_\CC \,Q_h$ is, of course, Arnold's Milnor number $\mu_\p (h)$ of $h$ considered as a function germ on $(\CC^{m+1},0)$ with the boundary $x_{3,4}=0$ or $x_{22}=0$ \cite{Ab}. In the quasi-homogeneous case $\mu_\p (h)$ coincides with $\tau_\p (h),$ the Tjurina number of the function with respect to the contact version of the boundary equivalence. The latter appears in the matrix context on its own if we consider the $\GL$-equivalence within the families  (\ref{Esymb}), (\ref{Esqb}) and (\ref{Eskb}): it  has already been noticed in \cite{B,BT} that for these matrix singularities

\centerline{
$ \tau_{\GL,\Mat_n} (M) = \tau_\p (h) \,,$
}

\medskip
\noindent
without any quasi-homogeneity requirements.
Indeed, the $\GL$-equivalence splits the last generator of the ideal in (\ref{Ebound}) into its two summands. 
}
\end{rem}

\begin{problem} {\em
Give an example of a function $h \in \O_{m+1}$ for which ${\rm dim}_\CC \,Q_h \ne \mu_\p (h).$
}
\end{problem}

\section{Geometry of bifurcation diagrams} \label{Sgbd}

\subsection{Bifurcation diagram of a matrix family} \label{Sbd}
Let $M: (\CC^s,0) \to \Mat_n$ be a germ of a matrix family with a finite $\SL$ Tjurina number $\tau,$ and
$(\LL,0) \simeq (\CC^\tau,0)$ the germ of the base of its $\SL$-miniversal deformation $\{M_\l\}_{\l \in (\LL,0)}$, $M_0=M.$  

The discriminant $\Delta$ in the target matrix space $\Mat_n$ is stratified by the sets of matrices of fixed coranks $c$: $\De = \cup_{c=1}^n \De_c.$ For a generic $\l \in \LL,$ the map $M_\l$ is transversal to this stratification.

\begin{definition} {\em
The germ $(\S,0) \subset (\LL,0)$ of the set of all those values of the deformation parameter $\l \in \LL$ for which the maps $M_\l$ are not transversal to the stratification of the discriminant $\De$ will be called}
the bifurcation diagram of the matrix family $M.$ 
\end{definition}

\begin{definition} {\em 
The {\em monodromy group\/} of a matrix family $M$ is the image of the representation of $\pi_1(\LL \setminus \S)$ on the middle homology of the singular Milnor fibre of $M.$
}
\end{definition} 

A generic choice of $\l \in \S$ gives us the mapping $M_\l$ with one most generic non-transversali\-ty to one of the strata $\De_c,$ that is, with one singularity of the $\SL$-stable class $A_1^{\Mat_c}.$ The Tjurina number of each of these singularities is 1, and therefore the diagram $\S$ is a hypersurface which splits into components according to the coranks: 
$$\S = \cup_{c=1}^n \S_c.$$

When a point $\l \notin \S$ approaches a generic point of $\l_* \in \S_c,$ a corank $c$ vanishing cycle in the singular Milnor fibre $ \De(M_\l)=(\det \circ M_\l)^{-1}(0)$ of the matrix family $M=M_0$ contracts to a point in $\De(M_{\l_*})$. (At the same moment a Morse critical point of the determinantal function lands on its zero level, which reduces by 1 the total Milnor number of critical points of the determinantal function outside its zero level.) Therefore, the diagram $\S$ may also be characterised as the set of all points $\l \in \LL$ for which the fibre $(\det \circ M_\l)^{-1}(0)$ is homotopic to a wedge of spheres of middle dimension whose number is smaller than the singular Milnor number $\mu_\De(M)$ of the family $M$. 
 
\begin{example} \label{EXbd} {\em
Let $M$ be a symmetric matrix family of the form (\ref{Esymb}) whose $\SL$ Tjurina number is $\tau.$ To obtain its $\SL$-miniversal deformation $\{ M_\l \},$ we replace the function $h(x_{22},z)$ in  (\ref{Esymb}) by a deformation $H(x_{22},z,\l) = H_\l(x_{22},z)$ of $h=H_0$ such that its initial velocities $\p H/ \p \l_1|_{\l=0},\dots,\p H/ \p \l_\tau|_{\l=0}$ form a basis of the symmetric version of the local algebra $Q_h$ from (\ref{Ebound}). Similar to (\ref{Ecrbound}), critical points of the function $\det \circ M_\l,$ $\l \notin \S,$ outside its zero level are critical points of the product
\begin{equation} 
\left( H_\l(x_{22},z)/(n-1) \right)^{n-1} \cdot x_{22}
\label{Efactors}
\end{equation}
outside its own zeros, that is, are given by the conditions
\begin{equation} 
{\rm grad}_z H_\l = 0 \quad {\rm and} \quad H_\l + (n-1) x_{22}\p H_\l/ \p x_{22} = 0\,.
\label{Ecrit}
\end{equation}
For $\l \in \S,$ some of these critical points land on the zero level of the product.
There are two options for such degenerations:
\begin{itemize}
\item[] either $x_{22}=0,$ and hence ${\rm grad}_z H_\l|_{x_{22}=0} = 0$ and  $H_\l|_{x_{22}=0} = 0,$
\item[] or $H_\l=0$ and ${\rm grad}_{x_{22},z} H_\l = 0.$
 \end{itemize}
For each option, the number of vanishing factors in (\ref{Efactors}) is the corank of the vanishing cycle:
the diagram components are respectively $\S_{n}$ and $\S_{n-1}.$

We notice that the diagrams $\S$ in this example are exactly what they should be if $H$ is considered as a deformation of the function germ $h$ on $(\CC^{m+1},0)$ with the boundary $x_{22}=0:$ 
the component $\S_n$ corresponds to the restriction of $H_\l$ to the boundary having critical value 0, while $\S_{n-1}$ corresponds to $H_\l$ itself having critical value 0 on $\CC^{m+1}.$
}
\end{example}

\begin{rem} \label{Rbound} {\em
For simple quasi-homogeneous matrix families $X_\tau^{\Mat_n}$ of Theorem \ref{Tsimple} 
$$
\tau = \tau_{\GL,\Mat_n} = \tau_{\SL,\Mat_n} = \mu_\De\,.
$$
$\GL$- or equivalently $\SL$-miniversal deformations of such families may be obtained by replacing the functions $h$ in (\ref{Esymb}--\ref{Eskb}) by their $\R_\p$-miniversal deformations.}
 \end{rem}

\subsection{The $k(\pi,1)$ theorem} \label{Skpi1}
We now enlarge the matrix diagrams $\S$ and include bifurcations of critical values outside zeros of the determinantal and Pfaffian functions.

\begin{definition} \label{Dfbd} 
The full bifurcation diagram of a matrix family germ $M: (\CC^s,0) \to \Mat_n$ {\em is the germ of the subset $(\Theta,0)$ in the base $(\LL,0)$ of an $\SL$-miniversal deformation $\{M_\l\}$ of $M$ for which the function $\det \circ M_\l$ or $\Pf \circ M_\l$ has locally fewer than $\mu_\De(M)$ distinct non-zero critical values.
}
\end{definition}

Besides the diagram $\S,$ the full diagram $\Theta$ contains the {\em Maxwell stratum\/} when the function $\det \circ M_\l$ or $\Pf \circ M_\l$ has coinciding non-zero values at different critical points, and the {\em caustic\/} corresponding to non-Morse critical points outside $\De (M_\l).$

\medskip
Let  $\Pi_{d} \simeq \CC^{d}$ be the space of monic degree $d$ polynomials in one variable. 
For a matrix family $M,$ we have a Lyashko-Looijenga type map from $\LL \setminus \Theta$ to $\Pi_{\mu_\De}$ which sends a point $\l$ to the polynomial whose roots are all non-zero critical values of the function $\det \circ M_\l$ or $\Pf \circ M_\l.$ By the continuity, it extends to the holomorphic map $\L: \LL \to \Pi_{\mu_\De}$ which sends $\Theta$ to the set  $\Xi \subset \Pi_{\mu_\De}$ of all polynomials with either multiple or zero roots. 

For matrix families  (\ref{Esymb}--\ref{Eskb}), the map $\L$ is a map between two $\tau$-dimensional spaces.

\begin{prop}\label{Pbih}
For matrix singularities (\ref{Esymb}--\ref{Eskb}), the map $\L: \LL \to \Pi_\tau$ is a local biholomorphism outside $\Theta.$ 
\end{prop}

{\em Proof.} We consider the symmetric case only. The two other settings are absolutely similar to it.

We will be working with the SL-miniversal deformations $\{ M_\l\}$ of matrices (\ref{Esymb}) introduced in Example \ref{EXbd}. We will be using the reduced expression (\ref{Efactors}) of $\det \circ M_\l$ for finding its critical points and values.

At a point $\l \notin \Theta,$ the values $c_1(\l), \dots, c_\tau(\l)$ of $\det \circ M_\l$ at its isolated critical points $a_1(\l), \dots, a_\tau(\l) \in \CC^{m+1}_{x_{22},z}$  are distinct. Therefore, the regularity of $\L$ is equivalent to the local regularity of the map $C:\l \mapsto (c_1, \dots, c_\tau).$ According to (\ref{Efactors}),
$$
c_i(\l) = \left. {x_{22}H^{n-1}
\over (n-1)^{n-1}}\right|_{(a_i(\l),\l)} \quad \Longrightarrow \quad
{\p c_i \over \p \l_j}(\l) = \left[ {x_{22}H^{n-2}
\over (n-1)^{n-2}} \cdot {\p H\over \p \l_j}
\right]_{(a_i(\l),\l)} ,
$$
since \,grad$_{x_{22},z} (x_{22}H^{n-1})$ vanishes at a point $(a_i(\l),\l)$ due to the criticality of the $a_i(\l).$ The fractional factor in the last expression does not depend on $j$ and is not zero since $\l \notin \Theta.$ Hence the Jacobi matrix of the map $C$ is invertible if and only if $\{\p H/ \p\l_1, \dots, \p H / \p \l_\tau\}$ (with $\l$ fixed) is a basis of the linear space of functions on the set of points $\{ a_1(\l) , \dots, a_\tau(\l) \}.$

To show that the $\p H / \p \l_j$ indeed form such a basis, consider the map germ $(\CC^{m+1+\tau},0) \to (\CC^{m+1+\tau},0)$:
$$
(x_{22}, z, \l) \mapsto (H + (n-1)x_{22} \p H/ \p x_{22}, \p H/ \p z_1, \dots, \p H/\p z_m,\l)
$$
(cf. \cite{Lya2}, section 2.4). Its local algebra is $Q_h,$ and the $\p H / \p \l_j$ represent a basis of $Q_h$ over $\CC.$ According to the Weierstrass-Malgrange preparation theorem \cite{Mal}, any holomorphic function germ in $x_{22}$ and $z$ may be written as a linear combination of the $\p H / \p \l_j$ with holomorphic coefficients:
$$
\sum_{j=1}^\tau \a_j \left(H + (n-1)x_{22} \p H/ \p x_{22}, \p H/ \p z_1, \dots, \p H/\p z_m,\l\right) \cdot \p H / \p \l_j\,\,.
$$
At the critical points $a_i(\l),$ the first $m+1$ arguments of the coefficients $\a_j$ vanish.
For a fixed $\l \notin \Theta,$ this yields a representation of a function on the set  
 $\{ a_1(\l) , \dots, a_\tau(\l) \}$ as a linear combination of the $\p H / \p \l_j$ with constant coefficients.
\hfill{$\Box$}

\bigskip
The classical Lyashko-Looijenga theorem \cite{L, Lya,Lya2} concerns the complement to the bifurcation diagram of a simple isolated boundary function singularity.
Proposition \ref{Pbih} is crucial for the following generalisation of that theorem to the simple matrix singularities $X_\tau^{\Mat_n}$ of Theorem \ref{Tsimple}.

\begin{theorem}\label{TLL} 
For a simple matrix singularity $X_\tau^{\Mat_n},$ the map $\L: \LL \to \Pi_{\tau}$ is a proper holomorphic map. Its restriction to the complement $\LL \setminus \Theta$ is a finite order unramified covering of the complement $\Pi_\tau \setminus \Xi.$
\end{theorem}

{\em Proof.} Once again, we are considering the symmetric case only.

A matrix family $X_\tau^{\Sym_n}$ is a particular case of the families (\ref{Esymb}), and we can use in its $\SL$-mniversal deformation described in Example \ref{EXbd} a quasi-homogeneous $\R_\p$-miniversal deformation $H$ of the function $h \in X_\tau.$ This gives us a quasi-homogeneous matrix deformation, with all the variables of positive weights.  The mapping $\L$ is therefore also quasi-homogeneous, with all its components of positive degrees. So, we consider $\L$ as a global map from $\LL = \CC^\tau$ to $\Pi_\tau = \CC^\tau.$ Its properness would follow from $\L^{-1} (0) = \{0\}.$

Assume $\l_* \in \L^{-1} (0),$ and consider a path in $\CC^\mu \setminus \Theta$ leading to $\l_*.$
For $\l \notin \Theta,$ all $\mu$ critical points $a_i(\l)$ of $\det \circ M_\l$ used in the construction of the map $\L$ are critical points of the function $x_{22} H_\l^{n-1}$ outside its zero level.
Approaching $\l_*,$ all these points land on the zero level of the function. Due to Example \ref{EXbd}, the latter means that all critical values of the function $H_{\l_*}$ considered as a function on $\CC^{m+1}$ with the boundary $x_{22} = 0$ must also be zero. The latter implies $\l_* = 0$ according to \cite{L,Lya,Lya2}. 

The Theorem follows now from Proposition \ref{Pbih}. \hfill{$\Box$}

\bigskip
The complement $\Pi_\tau \setminus \Xi$ is a $k(\pi,1)$-space for Brieskorn's generalised braid group $BrB_\tau$ associated with the Weyl group $B_\tau$ \cite{Br}. Therefore we have

\begin{corollary} \label{Ckpi1} 
For a simple matrix singularity $X_\tau^{\Mat_n},$ the complement $\LL \setminus \Theta$ to its full bifurcation diagram is a $k(\pi,1)$-space, where $\pi$ is a subgroup of a finite index in the braid group $BrB_\tau.$ In the symmetric and square matrix cases, the index is
$$
\left( (n-1) \C + \a \right)^\tau \tau ! / |X_\tau|\,.
$$
Here $\C$ and $|X_\tau|$ are the Coxeter number and the order of the Weyl group $X_\tau,$ and the values of $\alpha$
are as follows:
$$
\begin{array}{c||c|c|c|c}
X_\tau & A_\tau, D_\tau, E_\tau & B_\tau & C_\tau & F_4 \\  \vspace{-14pt} &&&& \\  \hline \vspace{-12pt} &&&& \\
\a & \C & 2 & 2\tau -2 & 6
\end{array} \,.
$$
In the skew-symmetric case, $n$ should be replaced in the formula by $k=n/2.$
\end{corollary}

The index has been calculated here as the degree of the quasi-homogeneous map $\L.$ In particular, formula (\ref{Efactors}) tells that $\alpha$ is the degree of the equation of the boundary when the degree of the function $X_\tau$ is set to be $\C.$  

\medskip
In Section \ref{Scv} we will extend the Theorem and its Corollary to the remaining simple matrix singularities of Theorem \ref{Tsimple}.

\section{Monodromy of corank 2 symmetric families} \label{S2x2}
Let us go back to the miniversal deformations (\ref{Evs}) of codimension 1 matrix singularities $A_1^{\Mat_n}.$
The quasi-homogeneous lift of the loop $\l=e^{2\pi i t}, t \in [0,1],$ from the base $\LL = \CC$ to the families of their singular Milnor fibres ends up with the map $(x,z) \mapsto (x,-z)$ of the domain $\CC^{N-1}_x \times \CC^m_z.$ Therefore, in all these cases the action of the only Picard-Lefschetz operator on the only vanishing cycle of the singular Milnor fibre is multiplication by $(-1)^m.$ Description of the monodromy representation of $\pi_1(\LL \setminus \S)$ on the homology of singular Milnor fibres of more complicated matrix families requires a good definition of an intersection form on the homology. In this section we are introducing an approach which allows to understand the whole monodromy group of
a corank 2  symmetric matrix singularity. 

\medskip
Coincidence of simple classifications of corank 2 symmetric matrix families and of functions on manifolds with boundary suggests introduction of an order 2 covering of the singular Milnor fibres. However, our covering will be different from the one used by Arnold in \cite{Ab} for boundary functions. 

Namely, since degenerate binary quadratic forms are squares of linear forms, we represent the cone $\De \subset \Sym_2$ as a quotient of a plane $\CC^2_{a,b}$ by the antipodal involution $\s: (a,b) \to 
(-a,-b).$ 
To a matrix family 
$$M: (\CC^s_u,0) \to (\Sym_2,0), \quad 
M(u)= \left(
\begin{array}{cc}
m_{11}(u) & m_{12}(u) \\
m_{12}(u) & m_{22}(u)
\end{array}
\right) 
$$
with a finite $\GL$ Tjurina number this associates a $\s$-invariant isolated complete intersection singularity 
\begin{equation}
\label{Ecov}
\wt M:  \ \ \left(
\begin{array}{cc}
m_{11}(u) & m_{12}(u) \\
m_{12}(u) & m_{22}(u)
\end{array}
\right) =
 \left(
\begin{array}{cc}
a^2 & ab \\
ab & b^2
\end{array}
\right) \quad {\rm in} \quad \CC^{s+2}_{u,a,b}\,.
\end{equation}
The isolatedness of the singularity of $\wt M$ follows from the transversality of the map $M$ to the stratified variety $\De$ on $\CC^s \setminus \{0\}.$ The latter holds due to $\tau_{\GL,\Sym_2}(M) < \infty$ (see, for example, \cite{DPI, DPII}).

Doing the same for a whole $\SL$-versal deformation of $M,$ we lift the singular Milnor fibre $V$ of $M$ 
to the smooth $\s$-invariant Milnor fibre $\wt V$ of $\wt M.$ The 2-covering $\wt V \to V$ is ramified  over the singularities of $V.$

This way the vanishing of a corank 2 cycle from Section \ref{Smvc}
$$
 \left(
\begin{array}{cc}
-x - z_1^2 - \dots - z_{s-2}^2 + \l & y \\
y & x
\end{array}
\right) 
$$
lifts to the standard vanishing 
$
a^2 + b^2 + z_1^2 + \dots + z_{s-2}^2 = \l 
$
of a Morse cycle in $\CC^s.$  We call this $\s$-invariant Morse cycle a {\em short cycle.}

On the other hand, the lift of a corank 1 vanishing cycle consists of two disjoint Morse cycles $e$ and $\s_* (e).$ We call their sum $e + \s_*(e)$ a {\em long cycle.}

The self-intersections of short and long cycles are the same as in Arnold's boundary function situation: zeros if $s$ is even,  $2$ and respectively $4$ if $s \equiv 1\mod4,$ and $-2$ and respectively $-4$ if $s \equiv 3\mod4.$

\medskip
\begin{definition} \label{Dif} {\em  (cf. \cite{Ab})
The {\em intersection form} of a corank 2 matrix family $M: (\CC^s,0) \to (\Sym_2,0)$ is a pair consisting of
\begin{itemize}
\item  the $\s$-invariant part $H^+$ of $H_{s-1}(\wt V; \ZZ),$ and
\item the restriction to $H^+$ of the intersection form of $\wt V.$
\end{itemize}}
\end{definition}

The Picard-Lefschetz operators acting on $H^+$ and corresponding to short and long cycles $e$ are given by the same formulas as in Arnold's theory of functions on manifolds with boundary:
\begin{equation}
\label{Eplo}
c \mapsto c + (-1)^{s(s+1)/2} (c \circ e) e \qquad {\rm and \  respectively} \qquad
c \mapsto c + (-1)^{s(s+1)/2} (c \circ e) e/2\,,
\end{equation}
where $(c \circ e)$ is the intersection number.

\begin{example} \label{XA2} {\em
Consider 
a deformed $A_2^{\Sym_2}$ family
$M'(x,y,z) =\left(
\begin{array}{cc}
-x + z^3 - z & y \\
y & x
\end{array}
\right)
$ in three variables.
The surface $V = \De(M')=\{x(-x+z^3 -z) - y^2 = 0\}$ is homotopic to the wedge of two corank 2 vanishing cycles: its own real parts $e_-$ within $-1 \le z \le 0$ and $e_+$ within $0 \le z \le 1.$ Their only common point is the origin, a singular point of $V.$ 

The covering (\ref{Ecov}) produces the surface $\wt V: z^3 - z = a^2 + b^2$ homotopic to a wedge of two 2-spheres: the real $\wt e_-$ and the purely imaginary $\wt e_+,$ which are 2-covers of respectively $e_-$ and $e_+.$ The short cycles $\wt e_-$ and $\wt e_+$ meet transversally in $\wt V$ at the origin. 
}
\end{example}

\begin{example}\label{Xblift} {\em 
For a matrix family $M(x,y,z) =\left(
\begin{array}{cc}
h(x,z) & y \\
y & x
\end{array}
\right), \ z \in \CC^{s-2},
$ the hypersurface $\De(M)$ lifts to the hypersurface $h (b^2,z) - a^2 = 0$ in $\CC^s_{a,b,z}.$
The latter is Arnold's lift $h(x,z)|_{x=b^2}$ of the boundary function $h$ followed by a one-variable stabilisation. The action of Arnold's boundary involution $b \to -b$ is extended here to the sign change of the stabilising variable (cf. \cite{Sy}). We have already noted that replacement of the function $h$ in $M$ by its $\R_\p$-versal deformation
provides a $\GL$-versal deformation of $M,$ with the bifurcation diagram of $h$ becoming the matrix bifurcation diagram $\S.$ Therefore, the monodromy group of the matrix family $M$ coincides in this case with the monodromy group of the one-variable stabilisation of the boundary function $h$.
}
\end{example}

The method we have introduced for understanding the monodromy of corank 2 symmetric matrix families  may be applied to higher corank singularities in the following way.

The singular Milnor fibre $V$ of a corank $n$ matrix family of any of the three types we have been considering in this paper has a filtration 
$$
V_n \subset V_{n-1} \subset \dots \subset V_2 \subset V_1 = V\,,
$$
where $V_c$ is 
the set of matrices of corank at least $c.$  The monodromy action of $\pi_1 (\LL \setminus \S)$ respects this filtration. The action of $\pi_1 (\LL \setminus \S)$ on $H_*(V_c)$ reduces to that of its quotient $\pi_1(\LL \setminus \cup_{i=c}^n \S_i).$

In the case of corank $n$ symmetric matrix families $M(u) = (m_{ij}(u))_{i,j=1}^n,$ $u \in \CC^s,$ our corank 2 approach generalises to the description of the monodromy on $V_{n-1}$ by introduction of an isolated complete intersection singularity
\begin{equation}
\label{E34}
\wt M:  \  (m_{ij}(u)) = (a_i a_j)  \quad {\rm in} \quad \CC^{s+n}_{u,a}
\end{equation}
symmetric under the involution $\s: a \to -a.$ The subset $V_{n-1}$ of the singular Milnor fibre of $M$ is 2-covered by the $\s$-invariant Milnor fibre $\wt V_{n-1}$ of $\wt M.$ The definition of the intersection form on $V_{n-1}$ follows Definition \ref{Dif} with the only difference that we should now take the part $H^{(-)^n}$ of the homology of $\wt V_{n-1}$ on which $\s$ acts as multiplication by $(-1)^n.$ Respectively, long cycles change to $e + (-1)^n \s_*(e).$   In the self-intersection numbers of short and long cycles and in the expressions (\ref{Eplo}) of the Picard-Lefschetz operators, the dimension $s$ should be replaced by $s+n -2.$ 

For symmetric matrix singularities of the form (\ref{Esymb}), this construction produces on $H_*(V_{n-1})$ either odd or even (depending on the dimension) versions of the monodromy groups  corresponding to the isolated boundary function singularities $h.$

\section{Corank 3 simple symmetric matrices \\ and simple odd functions} \label{Scork3}
Matrix families we have been considering so far are closely related to isolated functions singularities, at most boundary, and this relationship as a whole was spotted already during the simple classification carried out in \cite{B,BT}. We are now switching to the kind of `mysterious' part of the simple classification, the one from items (b) and (c) of Theorem \ref{Tsimple}. The approaches we are developing here for understanding these simple matrix singularities, are certainly applicable to the study of a much wider range of matrix families, especially of their bifurcations, topology and monodromy.  

\medskip
The aim of this section is to reveal a natural relation between three sets of objects:
\begin{itemize}
\item[(i)]
corank 3 simple symmetric matrix singularities in 4 variables, those from item (b) of Theorem \ref{Tsimple},
\item[(ii)]
simple singularities of odd functions on $\CC^2,$ 
\item[(iii)]
some of simple centrally symmetric complete intersection curves in $\CC^3.$ 
\end{itemize}
This relation was suggested by coincidence of the sets of degrees of basic invariants of certain Shephard-Todd groups and the sets of weights of parameters of miniversal deformations of such matrix families. Table \ref{T34} lists these quasi-homogeneous versal matrix deformations along with the corresponding odd functions $M'$ and curves $\wt M \subset \CC^3.$  

\begin{table} [htb]
\caption{Corank 3 simple symmetric matrix families in 4 variables} \vspace{-10pt}
$$\begin{array}{c|c|c|c}
\hline
\vspace{-10pt} &&& \\
M &  \SL{\rm -miniversal \  deformation \  of \  } M &  
\begin{array}{c} 
{\rm odd} \  M'\! : \\  \  \CC^2 \to \CC
\end{array} &
\begin{array}{c}
{\rm icis} \\ \wt M \subset \CC^3
\end{array}
\\
\vspace{-12pt} &&& \\
\hline \hline
\vspace{-10pt} && \\
\!\!\!\begin{array}{c}
{\rm I}_{k+1}, \\  k \ge 1 
\end{array} \!\!\! & 
\left( \begin{array}{ccc}
x      &     \l_k    & z \\
\l_k   & y + x^k + \l_{k-1}x^{k-1} + \dots + \l_1 x + \l_0 & w \\
 z     &      w     &  y
\end{array} \right) 
& 
\!\!\begin{array}{c}
D_{2k+2}/\ZZ_2 \\ \vspace{-1pt} \\
ac^2\! + \! a^{2k+1} \\ \vspace{-8pt} {}
\end{array}\!\!  &
\begin{array}{c}
 S_{2k+3} \\   \vspace{-8pt} \\
\!\!\!c^2\! +\!  2bc \! +\!  a^{2k}\! \! \!\\
ab
\end{array}
\\
\vspace{-11pt} &&& \\
\hline
\vspace{-10pt} && \\
{\rm II_4} 
&
\left( \begin{array}{ccc}
x      &  w^2 \! + \! \l_1 w \! + \! \l_0     & y \! + \! \l_3 w \! + \! \l_2 \\
w^2 \! + \! \l_1 w \! + \! \l_0 & y & z \\
y \! + \! \l_3 w \! + \! \l_2     &      z          & w
\end{array} \right) 
& 
\begin{array}{c}
E_8/\ZZ_2 
\\  \vspace{-1pt} \\
b^3 + c^5 \\ \vspace{-8pt} {}
\end{array} 
&
\begin{array}{c}
 U_9 \\   \vspace{-8pt} \\
b^2-ac \\
ab - c^4
\end{array}
\\
\vspace{-11pt} &&& \\
\hline
\vspace{-10pt} &&& \\
{\rm II_5} 
&
\!\!\!\left( \begin{array}{ccc} 
x      &\!\!\!\!\!\!\! \l_2 w^2 \! + \! \l_1 w \! + \! \l_0    & 
y\! + \! w^2 \! + \! \l_4 w \! + \! \l_3 
\!\!\!  \\ 
\!\! \l_2 w^2 \! + \! \l_1 w \! + \! \l_0 & y \  & z \!\!\!   \\
\!\!
y\!+\! w^2 \! + \! \l_4 w \! + \! \l_3    
&      z  \       & w \!\!\!
\end{array} \right) \!\!\!
 &
\begin{array}{c}
J_{10}/\ZZ_2  \\ \vspace{-1pt} \\
b^3 - bc^4  \\ \vspace{-8pt} {}
\end{array} 
& 
\begin{array}{c}
{U_{11}} \\  \vspace{-8pt} \\
\!\! b^2 \! - \! ac \! + \! c^4 \!\! \\
ab
\end{array}
\\
\vspace{-11pt} &&& \\
\hline
\vspace{-10pt} &&& \\
\vspace{-11pt}
{\rm II_6} 
&
\!\! \left( \begin{array}{ccc} \vspace{5pt}
 \!\!\!\!\!\! \! x   & \!\!\!\!\!\!\!\!\!\!\!\!\! \begin{array}{c}
w^3 \\
+ \l_2 w^2 \! + \! \l_1 w \! + \! \l_0
\end{array}     &\!\!\! \begin{array}{c}
y \   \\ + \l_5 w^2 \! + \! \l_4 w \! + \! \l_3 \end{array} \!\!\! \!\!\!
\\ \vspace{5pt}
\!\!\!\!\!\! \begin{array}{c}
w^3 \\ + \l_2 w^2 \! + \! \l_1 w \! + \! \l_0
\end{array}  & y \ \ \ \ \ \ \ \   & z \\
\!\!\!\!\!\! \begin{array}{c}
 y \ \\ + \l_5 w^2 \! + \! \l_4 w \! + \! \l_3 \end{array}      &      z \ \ \ \ \ \ \ \           & w
\end{array} \right) \!\!
&
\begin{array}{c}
E_{12}/\ZZ_2 \\ \vspace{-1pt} \\
b^3 + c^7 \\ \vspace{-8pt} {} 
\end{array} 
& 
\begin{array}{c}
{U_{13}} \\ \vspace{-8pt} \\
b^2 - ac \\
ab - c^6
\end{array}
\\
 &&& \\
\hline     
\end{array}
$$
\label{T34} 
\end{table}

To clarify the contents of the last two columns of the Table, let $\R_{odd,p} \subset \R$ be the subgroup of diffeomorphism germs of $(\CC^p,0)$ commuting with the central symmetry involution $\s_p$ on $(\CC^p,0).$ 
It is easy to check that the list of odd functions of item (ii) given in Table \ref{T34} is indeed a complete list of $\R_{odd,2}$-simple classes of such functions. The notations of the odd functions used in the Table are taken from \cite{GB,GMa,HJ}.

The simplicity of the curves in item (iii) is with respect to the action of the subgroup $\K_{odd,3} \subset \K$ in which the action of $\R$ on $\CC^3$ is restricted by that of $\R_{odd,3}.$  There are more $\K_{odd,3}$-simple curves in $\CC^3$ than we have in the Table. For example, $(a^2+b^4+c^4,bc)$ is $\K_{odd,3}$-simple (but not $\K$-simple).

The subscripts in the notation of the complete intersections $\wt M$ are their Milnor numbers.
The curve singularities $\wt M =  S_{2k+3}, U_9$ are 
from Giusti's list of simple curves in $\CC^3$ without any symmetry requirements \cite{Gi1,Gi2}.
They are the only centrally symmetric curves on that list.
The curves $\wt M = U_{11}, U_{13}$ become non-simple if the symmetry is removed. The way we denote them in the Table extends naturally Giusti's notation.

$\K_{odd,3}$-miniversal deformations of the curves $\wt M$ and $\R_{odd,2}$-miniversal deformations of the functions $M'$ may be obtained from respectively their monomial $\K$-miniversal and $\R$-miniversal deformations by omission of all odd (respectively even) monomials.

\subsection{Critical points and values}\label{Scv}
Our relationship between the matrix singularities and odd functions on $\CC^2$ will be based on comparison of their bifurcations. So, we introduce

\begin{definition} {\em 
Let $(\LL',0)$ be the base of an $\R_{odd,2}$-miniversal deformation of an odd holomorphic function germ $f$ on $(\CC^2,0).$
The {\em bifurcation diagram $ \S' \subset \LL'$} of $f$ is the germ of the set of all base points for which the zero level of the corresponding function is not smooth. The {\em full bifurcation diagram $ \Theta' \subset \LL'$} is the germ of the set of all base points for which the corresponding function is either not Morse or has a critical value zero. 
}
\end{definition}

All the singularities we are dealing with in this section are quasi-homogeneous as well as their miniversal deformations. Therefore, all our statements will be not about germs but about the whole complex coordinate spaces and algebraic varieties in them. 

\begin{prop} \label{Pfbd}
For each simple matrix family $M$ from Table \ref{T34} and the odd function $M'$ corresponding to it,
there is a biholomorphism between the bases of their miniversal deformations containing the full bifurcation diagrams:
$$ (\LL,\Theta) \simeq ( \LL', \Theta').
$$
\end{prop}

{\em Proof.} Our aim is to relate critical points of the function $\det \circ M_\l$ of particular members of a matrix versal deformation from Table \ref{T34} to critical points of the odd perturbations of the function $M'.$

In the case of the I$_{k+1}$ series, the first derivatives of $\det \circ M_\l$ with respect to $z$ and $w$ are twice the cofactors $C_{13}=C_{31}$ and $C_{23} = C_{32}$ of $M_\l.$
Remark \ref{Rblock} tells us that their vanishing at invertible critical matrices yields $m_{13}=m_{31}= m_{23} = m_{32} = 0,$ that is,
 $z=w=0.$ Thus at critical points off $\De$
$$
\det(M_\l)=
y \left| \begin{array}{cc}
x      &     \l_k   \\
\l_k   & y + P(x)
\end{array}
\right|,
$$
where $P$ is the polynomial in the central entry of the matrix in the table.
The vanishing of the $y$-derivative of $\det(M_\l)$ gives us $y=(\l_k^2-xP(x))/(2x)$ ($x=0$ would have implied $\l_k=0$ and the critical point landing on $\De$). Exclusion of $y$ reduces the determinant to a family of functions in just one variable:
$$\f (x;\l) = - \left(\l_k^2-xP(x)\right)^2/(4x).$$
The same function family comes from a miniversal deformation of the function singularity $D_{2k+2}$ in the class of odd functions on $\CC^2_{a,c}$:
$$
G(a,c;\l)=ac^2 + a P(a^2) + 2\l_k c, \quad {\rm where \ } P {\rm \  is \  the \  same \  as \  in \  I}_{k+1},
$$
when we start looking for critical points of its members outside their zero levels. Here exclusion of $c$ using $\p G/\p c =0$ reduces $G$ to
\begin{equation}\label{Esq1}
\psi (a;\l)  = \left(a^2 P (a^2) - \l_k^2\right)/a, \qquad {\rm that \  is, \  \ } \f(a^2;\l) = -\psi(a;\l)^2/4.
\end{equation}

The situation with the II$_\tau$ families is very similar, but more straightforward. Namely, for all of them
$\p (\det \circ M_\l)/\p x = C_{11}$ and --- in the spirit of Section \ref{S2x2} --- its vanishing  at a critical matrix allows us to set at the lower right corner of the matrix deformation
$$\left( \begin{array}{cc}
y & z \\
z  &  w
\end{array}
\right) :=
\left( \begin{array}{cc}
b^2 & bc \\
bc  &  c^2
\end{array}
\right),$$
where the $bc$-plane is equipped with the central symmetry involution $\s_2.$ This substitution yields
\begin{equation}\label{Esq2}
\det\circ M_\l = -\left(b^3 + bq(c^2) - cp(c^2)\right)^2,
\end{equation}
where $p$ and $q$ are the polynomials in $w$ in the matrix 
entries $m_{12}$ and $m_{13}$ in Table \ref{T34} including the constant $\l$-terms, but not $y.$ The expression in the brackets here is an $\R_{odd,2}$-miniversal deformation of the functions $E_8,$ $J_{10}$ and $E_{12}$ respectively, that is, exactly those mentioned in the third column of the table. 

The outcome of these calculations is as follows.

Each matrix singularity $M$ in the table is indexed by its $\SL$ Tjurina number $\tau.$ 
The ordinary Milnor number of each odd function $M'$ there is $2\tau.$ 
According to (\ref{Esq1}) and  (\ref{Esq2}), we can use the same space $\LL = \LL' = \CC^\tau$ for the bases of miniversal deformations of both $M$ and $M'.$ 
With such a choice:
\begin{itemize}
\item[(i)]
A generic member $M'_\l$ of the $\R_{odd,2}$-miniversal deformation of $M'$
has $\tau$ distinct pairs $(c_i,-c_i)$ of non-zero critical values of opposite signs taken at centrally symmetric pairs of critical points. 
Each pair corresponds to exactly one non-zero critical value of the determinant of the generic member $M_\l$ of the matrix deformation. All critical points in this case are Morse. 
\item[(ii)]
If $\l \in \Theta',$ then the odd function $M'_\l$ has coinciding non-zero critical values or degenerate critical points off its zero level if and only if the function $\det \circ M_\l$ has the same.
\item [(iii)]
When a generic $\l \notin \Theta'$ approaches a generic point $\l_*$ of $\S',$ 
a pair of centrally symmetric Morse critical points of the function $M'_\l$ ends up as a pair of distinct Morse points on $M'_{\l_*}=0.$ Respectively, one Morse critical point of $\det \circ M_\l$ lands on
$\det \circ M_{\l_*} = 0.$
\end{itemize} 
Therefore, we have $\Theta' = \Theta,$ and in particular $\S' = \S.$ \hfill{$\Box$}

\medskip
\begin{corollary}\label{Cmt}
For all simple matrix singularities of Theorem \ref{Tsimple}, 
\quad $\mu_\De = \tau_{SL, \Mat_n}. $
\end{corollary}

{\em Proof.} Singularities in parts (a) and (b) of Theorem \ref{Tsimple} are covered by respectively Theorem \ref{Tb} and item (i) by the end of the last proof. 

The only remaining case is  that of the simple seven-variable families $\wh M$ of $3\times 3$ square matrices. According to Theorem \ref{Tsimple}, each of them has the form
$\wh M = M + U,$ where $M$ is a simple four-variable symmetric family of Table 1, and the additional three variables are the entries of the skew-symmetric matrix 
$
U=
\left( \begin{array}{ccc}
0 & u_{12} & u_{13} \\
-u_{12} & 0 & u_{23}  \\
-u_{13} & -u_{23} & 0
\end{array}
\right).
$
It has been noted that the $\GL$ (same as $\SL$ in this case) Tjurina numbers of $\wh M$ and $M$ coincide.
Moreover, it is not so difficult to check (cf. \cite{BT}) that for an $\SL$ miniversal deformation of $\wh M$ we can take
$\{\wh M_\l = M_\l + U\}_{\l \in \LL},$ where $\{M_\l\}_{\l \in \LL}$ is an $\SL$ miniversal deformation of $M.$ Then 
$$
\p (\det \circ \wh M_\l)/ \p u_{ij} = C_{ij} - C_{ji}\,\,,
$$ 
the difference of the cofactors of $\wh M_\l.$ Vanishing of these three derivatives at an invertible matrix from a family $\wh M_\l$ implies that the inverse of this matrix is symmetric, and hence the matrix itself is also symmetric, that is, $U= 0.$ Thus invertible critical matrices in any family $\wt M_\l$ are exactly those in the family $M_\l,$  and therefore
$\mu_\De (\wh M) = \mu_\De (M) = \tau_{\SL, Sym_3}(M)  = \tau_{\SL, Sq_3}(\wh M).$ \hfill{$\Box$}

\bigskip
Corollary \ref{Cmt} allows us to extend the Lyashko-Looijenga type Theorem \ref{TLL} to all simple matrix singularities of Theorem \ref{Tsimple}. The proof of the Corollary implies that the full bifurcation diagrams and the Lyashko-Looijenga maps of a simple symmetric matrix singularity $M$ and of its square matrix version $M+U$ are the same. Therefore we are formulating our statement just for the symmetric matrices. Keeping the notation used in Section \ref{Skpi1}, we have

\begin{prop}\label{Pcov}
For the simple symmetric matrix singularities I$_\tau$ and II$_\tau,$ the  map $\L: \LL \to \Pi_\tau$ is a proper holomorphic map. Its restriction to $\LL \setminus \Theta$ is a finite order unramified covering of the complement $\Pi_\tau \setminus \Xi.$
\end{prop} 

{\em Proof.} An equivalent way to formulate the Proposition would be to state it for the Lyashko-Looijenga maps of $\R_{odd,2}$-simple functions which use the squares  of the critical values. Our argument will follow exactly this interpretation. 

According to the proof of Proposition \ref{Pfbd}, the map $\L$ is now constructed from  the squares of the critical values of odd perturbations of the functions $M'$ from Table \ref{T34}. Therefore, the local biholomorphicity of $\L$ outside $\Theta$ follows now directly from the similar local biholomorphicity of the Lyashko-Looijenga map of an isolated function singularity (see \cite{L}).   

A justification of $\L^{-1}(0) = \{0\}$ is now based on Gabrielov's theorem from \cite{Gab} stating that the Dynkin diagram of an isolated function singularity is connected, which implies that the $\L^{-1}(0)$ is the $\mu$=constant stratum in $\LL' \simeq \LL$ (here $\mu$ is the ordinary Milnor number of $M'$). In the $\R_{odd,2}$-miniversal deformations of each of the functions from Table \ref{T34} this stratum is just the origin. \hfill{$\Box$}

\begin{corollary}
For any of the simple symmetric matrix singularities I$_{\tau}$ and II$_\tau,$ the complement $\LL \setminus \Theta$ to its full bifurcation diagram is a $k(\pi,1)$-space, where $\pi$ is a subgroup of a finite index in the braid group $BrB_\tau.$ The indices are as follows:
$$
\begin{array}{c||c|c|c|c}
M & {\rm I}_{\tau} & {\rm II}_4 & {\rm II}_5 & {\rm II}_6 \\  \vspace{-14pt} &&&& \\  \hline \vspace{-12pt} &&&& \\
index & 2(2\tau-1)^\tau  & 2 \cdot 15^3 & 12^5 & 70 \cdot 21^4
\end{array} \,.
$$
\end{corollary}

\subsection{Bifurcations of matrices and of centrally symmetric curves}\label{Smcsc}
We will now approach the simple $3 \times 3$ symmetric matrix families in 4 variables from the point of view indicated at the end of Section \ref{S2x2}, that is, looking at bifurcations of the singular locus of the singular Milnor fibre of a matrix family.

For each of the families $M$ of Table \ref{T34},  the singular locus $V_{2,\l}$ of its singular Milnor fibre $V_\l,$ $\l \in \LL \setminus \S,$ is a smooth curve.  Following (\ref{E34}), we lift the $\SL$-miniversal deformations we had in the Table to $\CC^{4+\tau+3}_{x,y,z,w;\l;a.b.c}$ by equating them
$$
{\rm to \  }
\left( \begin{array}{ccc}
 a^2      &     ab      &  a(c+b) \\
  ab        &      b^2  &  b(c+b) \\
  a(c+b)  &    b(c+b)   &  (c+b)^2
\end{array} \right) 
\ {\rm for \  I}_\tau \quad {\rm and} \quad {\rm to \  }
\left( \begin{array}{ccc}
 a^2      &     ab      &  ac \\
  ab        &      b^2  &  bc \\
  ac        &      bc     &  c^2
\end{array} \right) \  {\rm for \  II_\tau}.
$$
The linear part of each of our map germs $M:(\CC^4,0) \to (\Sym_3,0)$ has rank 4, and therefore the lift yields deformations $\{ \wt V_{2,\l}\}$ of complete intersection curves $\wt M$ embeddable into $\CC^3_{a,b,c}$ and symmetric under the involution $\s_3: (a,b,c) \to (-a,-b,-c).$ The curves $\wt M$ form the last column of Table \ref{T34}. The lifted $\SL$-miniversal deformations $\{\wt V_{2,\l}\}$  turn out to be miniversal in the class of $\s_3$-symmetric complete intersections. In the notation used in the previous subsection, the deformations are 
\begin{equation}\label{Evd3}
\left( \begin{array}{c}
2bc + c^2 + P(a^2) \\
 ab - \l_k
\end{array} \right) 
\ {\rm for \  I}_{k+1} \qquad {\rm and} \qquad 
\left( \begin{array}{ccc}
 -ac + b^2 + q(c^2) \\
ab - p(c^2)
\end{array} \right) \  {\rm for \  II_\tau}.
\end{equation}

The subscripts in the Table notation of the complete intersections $\wt M$ are their Milnor numbers which all turn out to be equal to $2\tau + 1.$ The map $\wt V_{2,\l} \to V_{2,\l},$ $\l \notin \S,$ is an unramified 2-covering, and therefore, in all table cases a non-singular curve $V_{2,\l}$ is homotopic to a wedge of $\tau + 1$ circles. Notice that the whole singular Milnor fibre of each matrix family $M$ is homotopic to a wedge of just $\tau = \mu_\De$ copies of $S^3.$

\begin{definition} {\em 
Let $(\wt \LL,0)$ be the base of a $\K_{odd,3}$-miniversal deformation of a centrally symmetric one-dimensional icis germ $\wt M \subset (\CC^3,0).$ The {\em bifurcation diagram $\wt \S \subset \wt\LL$} of $\wt M$ is the germ of the set of all base points corresponding to singular curves. 
}
\end{definition}

Identifying the bases of the miniversal deformations of the matrix family $M$ and of the $\s_3$-symmetric curve $\wt M,$ we have $\wt \S (\wt M) \subseteq \S (M).$

\medskip
The deformations (\ref{Evd3}) bring us immediately to the $\R_{odd,2}$-miniversal deformations of the odd functions $M'$ 
which should be now considered as $\K_{odd,2}$-miniversal deformations. 

For example, in the II$_{\tau}$ case, equations (\ref{Evd3}) of the $\wt V_{2,\l}$ may be considered as a system of two linear equations in $a$ defining the unique value of $a,$ hence the determinant $b^3 + bq(c^2) - cp(c^2)$ of this system must be zero. The last expression is exactly the odd miniversal deformation of the function $M'$ we had in (\ref{Esq2}). There is only one point, the origin, on the determinantal curve $V'_\l = \{M'_\l=0\} \subset \CC^2_{b,c}$ for which we do not have any value of $a$ unless $p(0)=q(0)=0$, that is, when $\l$ is outside a codimension 2 plane -- denote it $\Upsilon$ -- in $\LL.$ Thus the projection $\CC^3_{a,b,c} \to \CC^2_{b,c}$ establishes a diffeomorphism $\wt V_{2,\l} \simeq V_\l' \setminus \{0\}$ for all $\l \in \LL \setminus \Upsilon.$ Under this projection, the central symmetry of $\CC^3$ becomes the central symmetry of $\CC^2.$ A reflection of the diffeomorphism between the curves is the relation between the Milnor numbers: $\mu(M') = \mu(\wt M) - 1 = 2\tau.$ 

The same happens in the I$_{k+1}$ case for the projection to the $ac$-plane. Therefore, for all pairs of singularities $M', \wt M$ from Table \ref{T34}, we have the coincidence of their bifurcation diagrams. 
Taking into account the last line of the proof of Proposition \ref{Pfbd}, we have

\begin{prop} \label{Pdiag}
For each triplet $M, M', \wt M$ of simple singularities from Table \ref{T34}, the pairs consisting of the base of a corresponding miniversal deformation and the bifurcation diagram of a singularity are biholomorphic:
 $$(\LL, \S) \simeq (\LL', \S') \simeq (\wt \LL, \wt \S).$$ 
\end{prop}

In fact, a fourth member may by added to each triplet in the Proposition so that its claim still stays valid, this time for all four bifurcation diagrams. This additional member is the simple $3 \times 3$ square matrix family $\wh M = M+U$ in seven variables.  

\medskip
The relation in the Proposition shows in particular that the diagram $\S(M)$ has just the corank 2 component $\S_2.$ Indeed, $\S_3(M) = \emptyset$  due to the dimension of the domain.  On the other hand, according to the previous subsection, all degenerations of the singular Milnor fibres $V_\l$ of $M$ are reflected in degenerations of the curves $M'_\l = 0,$ which in their turn correspond to degenerations of the curves $\wt V_{2,\l}$ and hence to degenerations of the singular loci $V_{2,\l} \subset V_\l.$ Hence $\S_1(M)$ is also empty.  

\begin{prop}
For all singularities of Table \ref{T34}, the multiplicities of their bifurcation diagrams are $\tau +1.$ 
\end{prop}

{\em Proof.} According to \cite{GB}, the bifurcation diagrams $\S'$ of the odd functions $D_{2k+2}/\ZZ_2$ and $E_8/\ZZ_2$ are the discriminants of the Shephard-Todd groups $G(4,2,k+1)$ and $G_{31}.$ Hence the multiplicities of $\S'$ in these cases are $\tau +1$ \cite{ST}.

The diagram $\S({\rm I}_{k+1}) \simeq \S'(D_{2k+2}/\ZZ_2)$ may be described as the set of all $\l \in \CC^{k+1}$ for which the polynomials
$$
xP(x) - \l_k^2 = x(x^k + \l_{k-1}x^{k-1} + \dots + \l_1 x + \l_0) - \l_k^2 
$$
have either multiple or zero roots \cite{ST}.

A common description of the diagrams $\S({\rm II}_\tau)$ is as follows (cf. \cite{GB}). 
Consider the centrally symmetric planar curves $b^3 + c p(c^2) - b q(c^2) = 0$ forming the miniversal family of the curve $M'=0.$ Blow up the plane at the origin  setting  $b=uc,$ take the strict transforms of the curves, and factorise by the symmetry action setting $c^2 = w.$ We end up with
\begin{equation}
u^3w + p(w) - u q(w) = 0\,.
\label{Epldef}
\end{equation}
In the three table cases this is a deformation of respectively $A_5,$ $D_6$ and $E_7.$ 
It extends to an $\R$-miniversal deformation by addition of the term $\l_\tau u^2.$

Therefore, the bifurcation diagrams of the families II${}_4,$ III${}_5$ and IV${}_6$ are non-generic sections of the discriminants of the Weyl groups $A_5,$ $D_6$ and $E_7$ by smooth hypersurfaces. The deformations (\ref{Epldef}) contain the constant and linear terms, which confirms that the multiplicities of the diagrams are indeed $\tau+1.$ \hfill{$\Box$}

\bigskip
I would like to finish with

\begin{conj}
Let $M: (\CC^s,0) \to \Mat_n$ be a germ of any of the three types of matrix families considered in this paper. Assume the number $s$ of its parameters is at least the codimension of the discriminant $\Delta$ in $\Mat_n,$ and the Tjurina number $\tau_{\SL,\Mat_n}(M)$ is finite. Then 

\centerline{$\mu_\Delta (M) = \tau_{\SL,\Mat_n}(M).$}
\end{conj}

We have shown that this is true for all known simple singularities and for the matrix versions of boundary function singularities. A similar amount of evidence existed behind the similar conjecture about the matrix families with few parameters which was proved successfully in \cite{GM}.

\ 

\noindent
{\em Department of Mathematical Sciences \\
University of Liverpool \\
Liverpool \  L69 7ZL  \\
UK}

\medskip
\noindent
{\tt goryunov@liverpool.ac.uk}


\begin{thebibliography}{99}

\bibitem{Ab} V.~I.~Arnold, {\em Critical points of functions on a manifold with boundary, the simple Lie groups $B_k, C_k, F_4$ and singularities of evolutes\/}, Russian Math. Surveys, {\bf 33} (1978), no. 5, 99--116. 

\bibitem{AGLV2} V.~I.~Arnold, V.~V.~Goryunov, O.~V.~Lyashko, and V.~A.~Vassiliev, {\em Singularities II. Classification and Applications\/}, Encyclopaedia of Mathematical Sciences, vol.39. Dynamical Systems VIII, Springer Verlag, Berlin a.o., 1993, v+256 pp.

\bibitem{AGV2} V.~I.~Arnold, S. M. Gusein-Zade, and A. N. Varchenko, {\em
Singularities of differentiable maps. Vol. II. Monodromy and asymptotics of integrals.}  Monographs in Mathematics {\bf 83}, Birkhäuser, 1988, viii+492 pp. 

\bibitem{Br}  E.  V.   Brieskorn,    
{\em  Sur  les  groupes  de  tresses  (d'apres  V.  I.  Arnold),} S\'eminaire Bourbaki, 24\`eme ann\'ee (1971/1972), Exp. No. 401, Lecture Notes in Math. {\bf 317} (1973), 21--44.

\bibitem{B} J.~W.~Bruce,
{\em  On families of symmetric matrices\/},  Mosc. Math. J. {\bf 3} (2003), no. 2, 335--360.

\bibitem{BGZ} J.~W.~Bruce, V.~V.~Goryunov, and V.~M.~Zakalyukin, 
{\em Sectional singularities and geometry of families of planar quadratic forms\/}, in {\em Trends in Singularities\/}, eds. A.~Libgober and M.~Tibar, Birkh\"auser, Basel, 2002, 83--97.

\bibitem{BT} J.~W.~Bruce and F.~Tari, 
{\em Families of square matrices\/}, Proc. London Math. Soc. {\bf 89} (2004), 738--762.

\bibitem{D} J.~Damon, {\em The unfolding and determinacy theorems for subgroups of $\A$ and $\K,$} Mem. Amer. Math. Soc {\bf 50} (1984), no. 306, x+88 pp.

\bibitem{Dh} J.~Damon, {\em Higher multiplicities and almost free divisors and complete intersections\/}, 
Mem. Amer. Math. Soc. {\bf 123} (1996), no. 589, x+113 pp.

\bibitem{Ds} J.~ Damon, 
{\em Schubert decomposition for Milnor fibers of the varieties of singular matrices\/}, J. Singul. {\bf 18} (2018), 358--396. 

\bibitem{DPI} J.~Damon, and B.~Pike, 
{\em Solvable groups, free divisors and nonisolated matrix singularities I: Towers of free divisors\/}, Ann. Inst. Fourier (Grenoble) {\bf 65} (2015), no. 3, 1251--1300.

\bibitem{DPII} J.~Damon, and B.~Pike, 
{\em Solvable groups, free divisors and nonisolated matrix singularities II: Vanishing topology\/}, Geometry \& Topology {\bf 18} (2014), 911--962.

\bibitem{Gab}
A.~M.~Gabrielov, 
{\em  Bifurcations, Dynkin diagrams and modality of isolated singularities\/}, Funkcional. Anal. i  Prilozhen. {\bf 8} (1974), no. 2, 7--12 (Russian). 
English translation: Functional Analysis and Its Applications {\bf 8} (1974), no. 2, 94--98.

\bibitem{Gi1} M. Giusti, 
{\em Classification des singularit\'es isol\'ees d'intersections compl\`etes simples\/}, C. R. Acad. Sci. Paris Sér. A-B {\bf 284} (1977), no. 3, A167--A170.

\bibitem{Gi2} M. Giusti, 
{\em Classification des singularit\'es isol\'ees simples d'intersections compl\`etes\/},  Singularities, Part 1 (Arcata, Calif., 1981), 457--494,
Proc. Sympos. Pure Math. {\bf 40}, Amer. Math. Soc., Providence, RI, 1983. 

\bibitem{GB} V.~Goryunov, and C.~Baines,
{\em Cyclic equivariant singularities of functions, and the unitary reflection groups $G(2m,2,n), G_9$ and $G_{31},$}  St. Petersburg Math. J. {\bf 11} (2000), no. 5, 761--774.

\bibitem{GMa} V.~Goryunov, and S.~H.~Man,
{\em The complex crystallographic groups and symmetries \linebreak of $J_{10},$} in:  Singularity theory and its applications, 
Adv. Stud. Pure Math. {\bf 43} (2006), 55--72.

\bibitem{GM} V.~Goryunov, and D.~Mond,
{\em Tjurina and Milnor numbers of matrix singularities\/}, J. London Math. Soc. {\bf 72} (2) (2005), 205--224.

\bibitem{GZ} V.~V.~Goryunov, and V.~M.~Zakalyukin, 
{\em Simple symmetric matrix singularities and the subgroups of Weyl groups $A\mu, D_\mu, E_\mu$},
Mosc. Math. J. {\bf 3} (2003), no. 2, 507--530.

\bibitem{HJ} J.~A.~Haddley,
 {\em Symmetries of unimodal singularities and complex hyperbolic reflection groups\/},
Ph.D. Thesis, University of Liverpool,  2011.

\bibitem{H} G.~Haslinger, 
{\em  Families of skew-symmetric matrices\/}, Ph.D. Thesis, University of Liverpool, 2001.

\bibitem{Le} D. T. L\^e, 
{\em Le concept de singularit\'e isol\'ee de fonction analytique,} Adv. Studies in Pure Math. {\bf 8} (1986), 215--227.

\bibitem{LG} D. T. L\^e, and G.-M. Greuel,
{\em Spitzen, Doppelpunkte und vertikale Tangenten in der Diskriminante verseller Deformationen von vollst\"andigen Durchschnitten\/}, Math. Ann. {\bf 222} (1976), 71--88.

\bibitem{L} E. Looijenga, 
{\em The complement of the bifurcation variety of a simple singularity,} Invent. Math. {\bf 23} (1974), 105--116.

\bibitem{Lya}   O. V. Lyashko, 
{\em The geometry of bifurcation diagrams,} Uspekhi Mat. Nauk {\bf 34} (1979), no. 3, 205--206 (Russian).

\bibitem{Lya2}  O. V. Lyashko, 
{\em Geometry of bifurcation diagrams,} Current problems in mathematics, Vol. 22, 94–129, Itogi Nauki i Tekhniki, Akad. Nauk SSSR, VINITI, Moscow, 1983 (Russian). English translation: Journal of Soviet Mathematics {\bf 27} (1984),  2736--2759.

\bibitem{Mal}
B. Malgrange, {\em Ideals of differentiable functions,} Oxford University Press, London, 1967, vii+106 pp.

\bibitem{ST} G. C. Shephard, and J. A. Todd,
{\em Finite unitary reflection groups\/}, Canad. J. Math. {\bf 6} (1954), 274--304.

\bibitem{S}  D.~Siersma,
{\em Vanishing cycles and special fibres\/}, Singularity Theory and its Applications, Warwick 1989, Part I,
Springer Lecture Notes in Mathematics {\bf 1462} (1991), 292--301.

\bibitem{Sy} P. Slodowy, 
{\em Simple singularities and simple algebraic groups,} Springer Lecture Notes in Mathematics {\bf 815} (1980), x+175 pp.

\end{thebibliography}
\end{document}